\documentclass[preprint]{elsarticle}

\usepackage{amsmath,amsthm,amssymb,amsfonts}
\usepackage{epsfig,epic}
\usepackage{paralist}
\usepackage{tabularx}
\usepackage{latexsym}
 \usepackage[utf8]{inputenc}
  \usepackage[colorlinks=true]{hyperref}
\hypersetup{urlcolor=blue, citecolor=red}

\newdefinition{ass}{Assumptions}[subsection]

\newtheorem{thm}{Theorem}[subsection]

\newdefinition{rmk}{Remark}
\newproof{pf}{Proof}
\newproof{pot}{Proof of Theorem \ref{thm2}}

\newdefinition{defin}{Definition}[subsection]

\begin{document}

\title{Hypocoercivity for a BGK model for gas mixtures}
\author[rvt1]{Liu Liu}
\ead{lliu@ices.utexas.edu}
\author[rvt2]{Marlies Pirner}
\ead{marlies.pirner@univie.ac.at}
\address[rvt1]{ICES, University of Texas at Austin, 201 E 24th St, Austin, TX, USA }
\address[rvt2]{University of Vienna, Oskar-Morgenstern-Platz 1, 1090 Wien}
\date{}

\begin{abstract}
We consider a kinetic model for a two component gas mixture without chemical reactions. Our goal is to study hypocoercivity for the linearized BGK model for gas mixtures in continuous phase space. By constructing an entropy functional, we can prove exponential relaxation to equilibrium with explicit rates. Our strategy is based on the entropy and spectral methods adapting Lyapunov's direct method as presented in \cite{achleitnerlinear} for the one species linearized BGK model. In comparison to the one species case, we start with two partial differential equations, one for each species. These equations are coupled due to interspecies interactions, which requires additional estimates on these interspecies terms.
\end{abstract}
\maketitle

\section{Introduction}

In this paper\textcolor{black}{,} we shall concern ourselves with a kinetic description of two gases. This is traditionally done via the Boltzmann equation for the two density distributions $f_1$ and $f_2$. Under certain assumptions the complicated interaction terms of the Boltzmann equation can be simplified by a so called BGK approximation (named after the physicists Bhatnagar-Gross-Krook \cite{BGK-Model}) consisting of a collision frequency multiplied by the deviation of the distributions from local Maxwellians. This approximation is constructed in a way such that it  has the same main properties of the Boltzmann equation namely conservation of the number of particles, momentum and energy. In addition\textcolor{black}{,} it has an H-theorem with an entropy inequality leading to an equilibrium which is a Maxwellian. For the BGK models, there are efficient numerical methods which are asymptotic preserving, 
meaning that the schemes remain efficient even approaching the hydrodynamic regime \cite{Puppo_2007, Jin_2010,Dimarco_2014, Bennoune_2008, Dimarco, Bernard_2015, Crestetto_2012}. The existence and uniqueness of solutions to the BGK equation for one species of gases in bounded domain in space was proven by Perthame and Pulvirenti in \cite{Perthame}.  
  
In this paper\textcolor{black}{,} we are interested in extension\textcolor{black}{s} of a BGK model to gas mixtures since in applications one often has to deal with mixtures instead of a single gas. \textcolor{black}{From the point of view of physicists, there are a lot of BGK models proposed in the literature concerning gas mixtures. Examples are the model of Gross and Krook in 1956 \cite{gross_krook1956}, the model of Hamel in 1965 \cite{hamel1965}, the model of Garzo, Santos and Brey in 1989 \cite{Garzo1989} and the model of Sofonea and Sekerka in 2001 \cite{Sofonea2001}. They all have one property in common.} Just like the Boltzmann equation for gas mixtures \textcolor{black}{that} contains a sum of collision terms on the right-hand side, \textcolor{black}{these kind\textcolor{black}{s} of models also have} a sum of collision terms in the relaxation operator. \textcolor{black}{In 2017, Klingenberg, Pirner and Puppo \cite{Pirner} proposed a kinetic model for gas mixtures which contains these often used models by physicists and engineers as special cases. Moreover, in \cite{Pirner} consistency of this model, like conservation properties (conservation of the number of particles of each species, conservation of total momentum and conservation of total energy), positivity and the H-Theorem, is proven. Since the models from physicists mentioned above are special cases of the model proposed in \cite{Pirner}, consistency of all these models is also proven. Another possible extension to gas mixtures was proposed by Andries, Aoki and Perthame in 2002 \cite{AndriesAokiPerthame2002}. In contrast to the other models it contains only one collision term on the right-hand side. Consistency like conservation properties (conservation of the number of particles of each species, conservation of total momentum and conservation of total energy), positivity and the H-Theorem is also proven there. Brull, Pavan and Schneider proved in \cite{Brull_2012} that the model \cite{AndriesAokiPerthame2002} can be derived by an entropy minimization problem. In recent works, there are efforts to extend this type of BGK model for gas mixtures to gas mixtures with chemical reactions, see for example the model of Bisi and C\'aceras \cite{Bisi}.  }

Once the existence and uniqueness of a steady state has been established, \textcolor{black}{one} can prove convergence to this steady state. 
It is more crucial and interesting to find quantitative estimates on the rates of convergence, 
which is known as hypocoercivity theory, see for example \cite{Villani, SGuo, CN, MB} for kinetic equations. 
There have been recent efforts extended to the study of kinetic equations 
with random inputs, including their mathematical properties such as regularity and long-time behavior in the random space, 
for example refer to \cite{HJ, LJ-UQ, LL1, LL2, Shu-J, DesPer}.
Although large-time behavior of the {\it monospecies} BGK equations were intensively studied in the literature, but \textcolor{black}{they} are unknown for 
{\it multispecies} BGK systems. The purpose of this paper is to study the large-time behavior of a linearized version of the BGK model for gas mixtures presented in \cite{Pirner}. We study hypocoercivity of this linearized model in one-dimensional phase space and construct an entropy functional to prove exponential relaxation to equilibrium with explicit rates. This paper is largely motivated by \cite{achleitnerlinear} and \cite{achleitner2017multi} for the one species BGK equation describing entropy and spectral methods adapting \textcolor{black}{Lyapunov's} direct method. In comparison to the one species case, we start with two partial differential equations, one for each species. These equations are coupled due to interspecies interactions. 
This requires additional estimates on these interspecies terms which will be proven in this work. These estimates are very different from the one species case because the interspecies terms from species $1$ and the interspecies terms from species $2$ have to be coupled and estimated together in an appropriate way. 
For example, we can not use conservation of momentum and energy of each species, since in gas mixtures we only have conservation of {\it total} momentum and energy. This requires a different way of estimating the total entropy that we carefully defined. 
 
The outline of the paper is as follows: In subsection 2.1, we will present the BGK model for two species developed in \cite{Pirner}. In subsection 2.2, we perform a linearization of this model assuming that the distribution functions are close to equilibrium. In subsection 2.3, we transform the system of partial differential equations to an infinite system of ordinary differential equations. In section 3, we define an appropriate entropy functional and develop additional estimates for hypocoercivity needed in the case of gas mixtures due to interspecies interactions. 
In section 4, we prove exponential relaxation with 
explicit estimates on the exponential convergence rate towards equilibrium. 

\section{Nonlinear and linearized BGK model for gas mixtures}
In this section we present the nonlinear and the linearized BGK model for gas mixtures, which is the main topic of this paper. 

\subsection{Nonlinear BGK model for gas mixtures}
We want to consider the BGK model for gas mixtures described in \cite{Pirner}. For the convenience of the reader, we want to briefly repeat it here. For more details, see \cite{Pirner}.
We consider the position space $\tilde{{\bf T}}^d:= \left( \frac{L}{2 \Pi} {\bf T}\right)^d$, the $d$-dimensional torus of side length $L$. Then, we consider the following BGK model for gas mixtures for two phase space densities $f_1(x,v,t), f_2(x,v,t); x \in \tilde{{\bf T}}^d, v \in \mathbb{R}^d$, one for each species, satisfying
\begin{align} \begin{split} \label{BGK}
\partial_t f_1 + v \cdot \nabla_x  f_1   &= \nu_{11} n_1 (M_1 - f_1) + \nu_{12} n_2 (M_{12}- f_1), \\[6pt]
\partial_t f_2 + v \cdot \nabla_x  f_2 &=\nu_{22} n_2 (M_2 - f_2) + \nu_{21} n_1 (M_{21}- f_2), \\
\end{split}
\end{align}
with the Maxwell distributions
\begin{align} 
\begin{split}
M_k &= \frac{n_k}{\sqrt{2 \pi \frac{T_k}{m_k}}^d }  \exp({- \frac{|v-u_k|^2}{2 \frac{T_k}{m_k}}})= \frac{n_k^{1+d/2}}{\sqrt{2 \pi \frac{P_k}{m_k}}^d }  \exp({- \frac{n_k|v-u_k|^2}{2 \frac{P_k}{m_k}}}), \quad k=1,2,
\\
M_{12} &= \frac{n_{1}}{\sqrt{2 \pi \frac{T_{12}}{m_1}}^d }  \exp({- \frac{|v-u_{12}|^2}{2 \frac{T_{12}}{m_1}}})= \frac{n_{1}^{1+d/2}}{\sqrt{2 \pi \frac{P_{12}}{m_1}}^d }  \exp({- \frac{n_1|v-u_{12}|^2}{2 \frac{P_{12}}{m_1}}}),
\\
M_{21} &= \frac{n_{2}}{\sqrt{2 \pi \frac{T_{21}}{m_2}}^d }  \exp({- \frac{|v-u_{21}|^2}{2 \frac{T_{21}}{m_2}}})= \frac{n_{2}^{1+d/2}}{\sqrt{2 \pi \frac{P_{21}}{m_2}}^d }  \exp({- \frac{n_2|v-u_{21}|^2}{2 \frac{P_{21}}{m_2}}}),
\end{split}
\label{BGKmix}
\end{align}
\textcolor{black}{where $m_k$, $k=1,2$ denotes the mass of a particle of species $k$.} The unknown quantities in these Maxwell distributions \eqref{BGKmix} will be explained in the equations \eqref{coll}--\eqref{temp}.
To be flexible in choosing the relationship between the collision frequencies, we now assume the relationship
\begin{equation} 
\nu_{12}=\varepsilon \nu_{21}, \quad 0 < \varepsilon \leq 1.
\label{coll}
\end{equation}
The restriction on $\varepsilon$ is without loss of generality. If $\varepsilon >1$, exchange the notation $1$ and $2$ and choose $\frac{1}{\varepsilon}$ as new $\varepsilon.$ For example, in the case of a plasma a common relationship found in the literature \cite{bellan2006} is given by $\varepsilon=\frac{m_1}{m_2}$.
The macroscopic quantities number density $n_k$, mean velocity $u_k$, temperature $T_k$, pressure $P_k$ are defined by
\begin{align}
\int f_k(v) \begin{pmatrix}
1 \\ v  \\ m_k |v-u_k|^2 
\end{pmatrix} 
dv =: \begin{pmatrix}
n_k \\ n_k u_k \\ d n_k T_k 
\end{pmatrix} , \quad P_k=n_k T_k \quad k=1,2.
\label{moments}
\end{align}
We define $P_{12}$ and $P_{21}$ by 
\begin{align} 
P_{12}=n_1 T_{12} \quad \text{and}  \quad P_{21}=n_2 T_{21}.
\label{P}
\end{align}
Then the remaining parameters  $ u_{12},$ $ u_{21},$ $ T_{12}$ and $T_{21}$ will be determined using conservation of  total momentum and energy, together with some symmetry considerations.
By choosing the densities of $M_{12}$ and $M_{21}$ equal to the denisities of the distribution functions $n_1$ and $n_2$, 
we have conservation of the number of particles, see Theorem 2.1 in \cite{Pirner}.
If we further \textcolor{black}{define} that $u_{12}$ is a linear combination of $u_1$ and $u_2$
 \begin{align}
u_{12}= \delta u_1 + (1- \delta) u_2, \quad \delta \in \mathbb{R},
\label{convexvel}
\end{align} then we have conservation of total momentum
provided that
\begin{align}
u_{21}=u_2 - \frac{m_1}{m_2} \varepsilon (1- \delta ) (u_2 - u_1),
\label{veloc}
\end{align}
see Theorem 2.2 in \cite{Pirner}.
If we further \textcolor{black}{define} that $T_{12}$ is of the following form
\begin{align}
\begin{split}
T_{12} &=  \alpha T_1 + ( 1 - \alpha) T_2 + \gamma |u_1 - u_2 | ^2,  \quad 0 \leq \alpha \leq 1, \gamma \geq 0 ,
\label{contemp}
\end{split}
\end{align}
then we have conservation of total energy
provided that
\begin{align}
\begin{split}
T_{21} =\left[ \frac{1}{d} \varepsilon m_1 (1- \delta) \left( \frac{m_1}{m_2} \varepsilon ( \delta - 1) + \delta +1 \right) - \varepsilon \gamma \right] |u_1 - u_2|^2 \\+ \varepsilon ( 1 - \alpha ) T_1 + ( 1- \varepsilon ( 1 - \alpha)) T_2,
\label{temp}
\end{split}
\end{align}
see Theorem 2.3 in \cite{Pirner}.
In order to ensure the positivity of all temperatures, we need to restrict $\delta$ and $\gamma$ to 
 \begin{align}
0 \leq \gamma  \leq \frac{m_1}{d} (1-\delta) \left[(1 + \frac{m_1}{m_2} \varepsilon ) \delta + 1 - \frac{m_1}{m_2} \varepsilon \right],
 \label{gamma}
 \end{align}
and
\begin{align}
 \frac{ \frac{m_1}{m_2}\varepsilon - 1}{1+\frac{m_1}{m_2}\varepsilon} \leq  \delta \leq 1,
\label{gammapos}
\end{align}
see Theorem 2.5 in \cite{Pirner}. 

Let $d\tilde{x}:= L^{-d} dx$ denote the normalized Lebesgue measure on ${\bf T}^d$. We consider normalized initial data
\begin{align}
\begin{split}
\int \int f_k^I(x,v) d\tilde{x} dv = 1, ~ \int \int v ( m_1 f_1^I + m_2 f_2^I ) d\tilde{x}dv = 0, \\[4pt] \int \int |v|^2 (m_1 f_1^I + m_2 f_2^I ) d\tilde{x} dv = n_{\infty,1} + n_{\infty,2}.
\end{split}
\label{consin}
\end{align}
We expect that equations \eqref{BGK}  have the unique space-homogeneous steady state 
$$ f_k^{\infty} (v) = \frac{n_{\infty,k}}{(2 \pi T_{\infty}/m_k)^{d/2}} \exp \left(- \frac{|v-u_{\infty}|^2}{2 T_{\infty}/m_k}\right), \quad k=1,2, $$ which are two \textcolor{black}{Maxwellian} with densities $n_{\infty,k}= \int \int f_k(x,v,0) dv d\tilde{x}$, equal mean velocity $u_{\infty}$ and equal temperature $T_{\infty}$.
By translating and scaling the coordinate system, we may assume $u_{\infty}=0$ and $T_{\infty}=1$ such that we obtain
\begin{align}
f_k^{\infty} (v) = \frac{n_{\infty,k}}{(2 \pi /m_k)^{d/2}} \exp \left(- \frac{|v|^2}{2 /m_k}\right), \quad k=1,2.
\label{equ}
\end{align}
\subsection{Linearized BGK model for gas mixtures}
In this section, we derive a linearized version of the BGK model for gas mixtures described in the previous section.

For this, we consider a solution $(f_1,f_2)$ to \eqref{BGK} which is close to the equilibrium $(f_1^{\infty}, f_2^{\infty})$ with 
\begin{align}
f_k(x,v,t) = f_k^{\infty}(v) + h_k(x,v,t).
\label{ansatz}
 \end{align}
Then, we have
\begin{align}
\begin{split}
 &n_k(x,t)= n_{\infty,k} + \sigma_k(x,t) \quad \text{with} \quad \sigma_k(x,t) = \int h_k(x,v,t) dv \\
&(n_k u_k)(x,t) = \int v f_k(x,v,t) dv = \mu_k(x,t) \quad \text{with} \quad \mu_k(x,t) = \int v h_k(x,v,t) dv \\ &P_k(x,t) = \frac{m_k}{d}  \int |v-u_k|^2 f_k(x,v,t) dv = n_{\infty,k} + \frac{1}{d}\left[\tau_k(x,t) - \frac{m_k |\mu_k(x,t)|^2}{n_{\infty,k}+ \sigma_k(x,t)}\right] \\ &\hspace{5cm}\quad \text{with} \quad \tau_k(x,t) = m_k \int |v|^2 h_k(x,v,t) dv.
\end{split}
\label{momh}
\end{align} 
The conservation of the normalization \eqref{consin} implies 
\begin{align}
\begin{split}
&\int \sigma_1(x,t) d\tilde{x} = \int \sigma_2(x,t) d\tilde{x} = 0, 
\\& \int (m_1 \mu_1(x,t) + m_2 \mu_2(x,t)) d\tilde{x} = 0, \,  \int (\tau_1(x,t) + \tau_2(x,t) ) d\tilde{x} = 0. 
\end{split}
\label{cons}
\end{align}
Now we derive the linearized version of the equations \eqref{BGK} by inserting the ansatz \eqref{ansatz} into \eqref{BGK}, 
then $h_1$ and $h_2$ satisfy
\begin{align} \begin{split} \label{BGKh}
\partial_t h_1 + v \cdot \nabla_x  h_1  &= \nu_{11} n_1 (M_1 - f_1^{\infty} - h_1) + \nu_{12} n_2 (M_{12}- f_1^{\infty} - h_1),
\\[8pt]
\partial_t h_2 + v \cdot \nabla_x h_2 &=\nu_{22} n_2 (M_2 - f_2^{\infty} -h_1) + \nu_{21} n_1 (M_{21}- f_2^{\infty} -h_2), \\
\end{split}
\end{align}

We want to linearize the model \eqref{BGKh} by performing a Taylor expansion of $M_1, M_2, M_{12}, M_{21}$ with respect to $\sigma_1, \sigma_2, \mu_1, \mu_2, \tau_1$ and $\tau_2$ around $0$ assuming that $\sigma_1, \sigma_2, \mu_1, \mu_2, \tau_1$ and $\tau_2$ are small. 
For the one-species terms, namely the first terms on the right-hand side of \eqref{BGKh}, we obtain
\begin{align*}
&\textstyle\quad M_k(x,v,t) -f_k^{\infty}(v) \\& \textstyle= \frac{(n_{\infty,k}+ \sigma_k(x,t))^{1+d/2}}{\sqrt{\frac{2 \pi}{m_k} ( n_{k, \infty} + \frac{1}{d}[\tau_k(x,t) -  \frac{m_k |\mu_k(x,t)|^2}{n_{\infty,k} + \sigma_k(x,t)}])}^d }\, \exp\left({- \frac{|v(n_{\infty,k} + \sigma_k(x,t))-\mu_k|^2}{ \frac{2}{m_k} ( n_{k, \infty} + \frac{1}{d}[\tau_k(x,t) -  \frac{m_k |\mu_k(x,t)|^2}{n_{\infty,k} + \sigma_k(x,t)}]) (n_{\infty,k} + \sigma_k(x,t))}}\right)
\\[4pt] & \textstyle\quad - \frac{n_{\infty,k}}{(2 \pi /m_k)^{d/2}} \exp (- \frac{|v|^2}{2 /m_k}), 
\\[4pt]& \textstyle \approx f_k^{\infty} (v) 
\left[(\frac{1+d/2}{n_{\infty,k}} - \frac{m_k |v|^2}{2 n_{\infty,k}}) \sigma_k(x,t) + \frac{m_k}{n_{\infty,k}} v \cdot \mu_k(x,t) + \frac{1}{n_{\infty,k}} (- \frac{1}{2} + \frac{m_k|v|^2}{2d}) \tau_k(x,t)\right]. 
\end{align*}
For the mixture part, we first recognize that $P_{12}$ defined in \eqref{P} and \eqref{contemp} and $\mu_{12}:= n_1 u_{12}$ can be written as 
$$ \mu_{12}= \delta \mu_1 + (1-\delta) \frac{n_1}{n_2} \mu_2, \quad P_{12}= \alpha P_1 + (1- \alpha) \frac{n_1}{n_2} P_2 + \gamma \frac{|\mu_1 - \frac{n_1}{n_2} \mu_2|^2}{n_1}.$$ If we insert this into the expression for $M_{12}$, we obtain
\begin{align*}
\scriptstyle
M_{12}(x,v,t) = \frac{n_{1}^{1+d/2}}{\sqrt{\frac{2 \pi}{m_1} (\alpha P_1 + (1-\alpha)\frac{n_1}{n_2} P_2 + \gamma \frac{|\mu_1 - \frac{n_1}{n_2} \mu_2|^2}{n_1}) }^d }\, \exp\left({- \frac{|v n_1- \delta \mu_1 - (1-\delta) \frac{n_1}{n_2} \mu_2|^2}{2 \frac{2 n_1}{m_1}(\alpha P_1 +(1-\alpha) \frac{n_1}{n_2} P_2 + \gamma \frac{|\mu_1 - \frac{n_1}{n_2} \mu_2|^2}{n_1})}}\right). 
\end{align*}
Next, we insert $n_k$ and $P_k$ given by \eqref{momh} and obtain 
\begin{align*}
&\quad M_{12}(x,v,t) \\&\scriptstyle= \frac{(n_{\infty,1} + \sigma_1)^{1+d/2}}{\sqrt{\frac{2 \pi}{m_1} (\alpha (n_{\infty,1} + \frac{1}{d}[\tau_1 - \frac{m_1 |\mu_1|^2}{n_{\infty,1}+ \sigma_1}]) + (1-\alpha)\frac{n_{\infty,1} + \sigma_1}{n_{\infty,2} + \sigma_2} (n_{\infty,2} + \frac{1}{d}[\tau_2 - \frac{m_2 |\mu_2|^2}{n_{\infty,2}+ \sigma_2}] + \gamma \frac{|\mu_1 - \frac{n_{\infty,1} + \sigma_1}{n_{\infty,2} + \sigma_2} \mu_2|^2}{n_{\infty,1} + \sigma_1}) }^d } 
\\& \scriptstyle \exp\bigg({- \frac{|v (n_{\infty,1} + \sigma_1)- \delta \mu_1 - (1-\delta) \frac{n_{\infty,1} + \sigma_1}{n_{\infty,2} + \sigma_2} \mu_2|^2}{ \frac{2 (n_{\infty,1} + \sigma_1)}{m_1}(\alpha (n_{\infty,1} + \frac{1}{d}[\tau_1 - \frac{m_1 |\mu_1|^2}{n_{\infty,1}+ \sigma_1}]) +(1-\alpha) \frac{n_{\infty,1} + \sigma_1}{n_{\infty,2} + \sigma_2} (n_{\infty,2} + \frac{1}{d}[\tau_2 - \frac{m_2 |\mu_2|^2}{n_{\infty,1}+ \sigma_2}]) + \gamma \frac{|\mu_1 - \frac{n_{\infty,1} + \sigma_1}{n_{\infty,2} + \sigma_2} \mu_2|^2}{n_{\infty,1} + \sigma_1})}}\bigg). 
\end{align*}
With this expression and by denoting the set $\mathcal D = \{\sigma_1=\sigma_2=\mu_1=\mu_2=\tau_1=\tau_2=0\}$, then 
one obtains the following derivatives
\begin{align*}
\partial_{\sigma_1} M_{12}|_{\mathcal D}&= \frac{1}{n_{\infty,1}}\left(1+ \frac{\alpha}{2} (d- m_1 |v|^2)\right),  \\
\partial_{\sigma_2} M_{12}|_{\mathcal D} &= \frac{1}{2} \frac{1}{n_{\infty,2}}(1-\alpha) (d- m_1 |v|^2), \\
\partial_{\mu_1} M_{12}|_{\mathcal D} &=  \frac{1}{n_{\infty,1}}\delta m_1 v,  \\
\partial_{\mu_2} M_{12}|_{\mathcal D} &=  \frac{1}{n_{\infty,2}}(1-\delta) m_1 v,  \\ 
\partial_{\tau_1} M_{12}|_{\mathcal D} &=  \frac{1}{2} \frac{1}{n_{\infty,1}}\alpha \left(\frac{1}{d} m_1 |v|^2 -1\right), \\
\partial_{\tau_1} M_{12}|_{\mathcal D} &=  \frac{1}{2} \frac{1}{n_{\infty,2}}(1-\alpha) \left(\frac{1}{d} m_1 |v|^2 -1\right). 
\end{align*}

Therefore, 
\begin{align*}
&M_{12} - f_1^{\infty} \approx f_1^{\infty} \Big[ \frac{1}{n_{\infty,1}}\left(1+ \frac{\alpha}{2} (d- m_1 |v|^2)\right) \sigma_1 + \frac{1}{2} \frac{1}{n_{\infty,2}}(1-\alpha) (d- m_1 |v|^2) \sigma_2 \\&\hspace{2.5cm}  + \frac{1}{n_{\infty,1}}\delta m_1 v \cdot \mu_1 + \frac{1}{n_{\infty,2}}(1-\delta) m_1 v \cdot \mu_2 + \frac{1}{2} \frac{1}{n_{\infty,1}}\alpha \left(\frac{1}{d} m_1 |v|^2 -1 \right) \tau_1 \\ &\hspace{2.5cm} +  \frac{1}{2} \frac{1}{n_{\infty,2}}(1-\alpha) \left(\frac{1}{d} m_1 |v|^2 -1 \right) \tau_2 \Big]. 
\end{align*}
Similarly, in the case of species $2$, we observe that $P_{21}$ defined in \eqref{P} and \eqref{temp} and $\mu_{21}= n_2 u_{21}$ can be written as
\begin{align*}
\mu_{21} = \tilde{\delta} \frac{n_2}{n_1} \mu_1 + (1- \tilde{\delta}) \mu_2, \quad P_{21}= \tilde{\alpha} \frac{n_2}{n_1} P_1 + (1- \tilde{\alpha} ) P_2 + \tilde{\gamma} \frac{|\frac{n_2}{n_1} \mu_1 - \mu_2|^2}{n_2}\,, 
\end{align*}
with $\tilde{\alpha}, \tilde{\delta}$ and $\tilde{\gamma}$ given by
\begin{align*}
\tilde{\alpha}= \varepsilon (1- \alpha), \quad \tilde{\delta}= \frac{m_1}{m_2} \varepsilon (1- \delta), \quad \tilde{\gamma}= \frac{1}{d} \varepsilon m_1 (1- \delta) ( \frac{m_1}{m_2} \varepsilon (\delta -1 ) + \delta +1) - \varepsilon \gamma\,. 
\end{align*}
Then, by inserting these expressions and the expressions of $n_1, n_2, P_1$ and $P_2$ given by \eqref{momh}, we obtain
\begin{align*}
&M_{21}(x,v,t) \\&\scriptstyle= \frac{(n_{\infty,2} + \sigma_2)^{1+d/2}}{\sqrt{\frac{2 \pi}{m_2} (\tilde{\alpha} \frac{n_{\infty,2} + \sigma_2}{n_{\infty,1} + \sigma_1}(n_{\infty,1} + \frac{1}{d}[\tau_1 - \frac{m_1 |\mu_1|^2}{n_{\infty,1}+ \sigma_1}]) + (1-\tilde{\alpha}) (n_{\infty,2} + \frac{1}{d}[\tau_2 - \frac{m_2 |\mu_2|^2}{n_{\infty,2}+ \sigma_2}] + \tilde{\gamma} \frac{|\frac{n_{\infty,2} + \sigma_2}{n_{\infty,1} + \sigma_1}\mu_1 -  \mu_2|^2}{n_{\infty,2} + \sigma_2}) }^d } 
\\& \scriptstyle \exp\bigg({- \frac{|v (n_{\infty,2} + \sigma_2)- \tilde{\delta} \frac{n_{\infty,2} + \sigma_2}{n_{\infty,1} + \sigma_1}\mu_1 - (1-\tilde{\delta})  \mu_2|^2}{ \frac{2 (n_{\infty,2} + \sigma_2)}{m_2}(\tilde{\alpha} \frac{n_{\infty,2} + \sigma_2}{n_{\infty,1} + \sigma_1} (n_{\infty,1} + \frac{1}{d}[\tau_1 - \frac{m_1 |\mu_1|^2}{n_{\infty,1}+ \sigma_1}]) +(1-\tilde{\alpha})  (n_{\infty,2} + \frac{1}{d}[\tau_2 - \frac{m_2 |\mu_2|^2}{n_{\infty,1}+ \sigma_2}]) + \tilde{\gamma} \frac{|\frac{n_{\infty,2} + \sigma_2}{n_{\infty,1} + \sigma_1}\mu_1 -  \mu_2|^2}{n_{\infty,2} + \sigma_2})}}\bigg). 
\end{align*}
From this we get the following derivatives
\begin{align*}
\partial_{\sigma_1} M_{21}|_{\mathcal D} &= \frac{1}{2} \frac{1}{n_{\infty,1}}\tilde{\alpha} (d- m_1 |v|^2) = \frac{1}{2} \frac{1}{n_{\infty,1}} \varepsilon(1-\tilde{\alpha}) (d- m_1 |v|^2), \\
\partial_{\sigma_2} M_{21}|_{\mathcal D} &= \frac{1}{n_{\infty,2}}\left(1+ \frac{1- \tilde{\alpha}}{2} \right) (d- m_1 |v|^2) 
= \frac{1}{n_{\infty,2}}\left(1+ \frac{1- \varepsilon (1- \alpha)}{2} \right) (d- m_1 |v|^2), \\
\partial_{\mu_1} M_{21}|_{\mathcal D} &=  \frac{1}{n_{\infty,1}}\tilde{\delta} m_2 v=  \frac{1}{n_{\infty,1}} \varepsilon (1- \delta) m_1 v, \\
\partial_{\mu_2} M_{21}|_{\mathcal D} &=  \frac{1}{n_{\infty,2}}(1-\tilde{\delta}) m_2 v =  \frac{1}{n_{\infty,2}}\left(1-\frac{m_1}{m_2}\varepsilon (1- \delta)\right) m_2 v, \\ 
\partial_{\tau_1} M_{21}|_{\mathcal D} &=  \frac{1}{2} \frac{1}{n_{\infty,1}}\tilde{\alpha} \left(\frac{1}{d} m_2 |v|^2 -1\right) =  \frac{1}{2} \frac{1}{n_{\infty,1}} \varepsilon (1- \alpha) \left(\frac{1}{d} m_2 |v|^2 -1\right), \\
\partial_{\tau_1} M_{21}|_{\mathcal D} &=  \frac{1}{2} \frac{1}{n_{\infty,2}}(1-\tilde{\alpha}) \left(\frac{1}{d} m_2 |v|^2 -1\right)
= \frac{1}{2} \frac{1}{n_{\infty,2}}(1-\varepsilon(1- \alpha)) \left(\frac{1}{d} m_2 |v|^2 -1\right). 
\end{align*}
Then, we get
\begin{align*}
&M_{21} - f_2^{\infty} \approx f_2^{\infty} \Big[\frac{1}{2} \frac{1}{n_{\infty,1}}\varepsilon(1-\alpha) (d- m_2 |v|^2) \sigma_1+ \frac{1}{n_{\infty,2}}\left(1+ \frac{1-\varepsilon(1-\alpha)}{2} (d- m_2 |v|^2)\right) \sigma_2   \\&  + \frac{1}{n_{\infty,1}}\varepsilon (1- \delta) m_1 v \cdot \mu_1 + \frac{1}{n_{\infty,2}}(1-\frac{m_1}{m_2} \varepsilon (1- \delta)) m_2 v \cdot \mu_2 + \frac{1}{2} \frac{1}{n_{\infty,1}}\varepsilon (1- \alpha) \left(\frac{1}{d} m_2 |v|^2 -1\right) \tau_1 \\ &   +  \frac{1}{2} \frac{1}{n_{\infty,2}}(1-\varepsilon (1- \alpha)) \left(\frac{1}{d} m_2 |v|^2 -1 \right) \tau_2 \Big]. 
\end{align*}

 To summarize, the linearized BGK equations are given by 
 {\footnotesize
\begin{align} \begin{split} \label{BGKh2}
&\quad\partial_t h_1 + v \cdot \nabla_x  h_1 \\ &= \nu_{11} (n_{\infty,1} + \sigma_1) \left(f_1^{\infty} (v) [(\frac{1+d/2}{n_{\infty,1}} - \frac{m_1 |v|^2}{2 n_{\infty,1}}) \sigma_1(x,t) + \frac{m_1}{n_{\infty,1}} v \cdot \mu_1(x,t) + \frac{1}{n_{\infty,1}} (- \frac{1}{2} + \frac{m_1|v|^2}{2d}) \tau_1(x,t)]- h_1\right) 
\\&\quad+ \nu_{12} (n_{\infty,2} + \sigma_2(x,t))\bigg(f_1^{\infty} [ \frac{1}{n_{\infty,1}}(1+ \frac{\alpha}{2} (d- m_1 |v|^2)) \sigma_1(x,t) + \frac{1}{2} \frac{1}{n_{\infty,2}}(1-\alpha) (d- m_1 |v|^2) \sigma_2(x,t) \\&\quad+ \frac{1}{n_{\infty,1}}\delta m_1 v \cdot \mu_1(x,t) + \frac{1}{n_{\infty,2}}(1-\delta) m_1 v \cdot \mu_2(x,t) + \frac{1}{2} \frac{1}{n_{\infty,1}}\alpha (\frac{1}{d} m_1 |v|^2 -1) \tau_1(x,t) \\&\quad+  \frac{1}{2} \frac{1}{n_{\infty,2}}(1-\alpha) (\frac{1}{d} m_1 |v|^2 -1) \tau_2(x,t)] - h_1\bigg), 
\\[10pt]
&\quad\partial_t h_2 + v \cdot \nabla_x h_2\\ &  =\nu_{22} (n_{\infty,2} + \sigma_2)\left(f_2^{\infty} (v) [(\frac{1+d/2}{n_{\infty,2}} - \frac{m_2 |v|^2}{2 n_{\infty,2}}) \sigma_2(x,t) + \frac{m_2}{n_{\infty,2}} v \cdot \mu_2(x,t) + \frac{1}{n_{\infty,2}} (- \frac{1}{2} + \frac{m_2|v|^2}{2d}) \tau_2(x,t)] -h_2\right)
\\&\quad  + \nu_{21} (n_{\infty,1} + \sigma_1(x,t))\bigg(f_2^{\infty}[\frac{1}{2} \frac{1}{n_{\infty,1}}\varepsilon(1-\alpha) (d- m_2 |v|^2) \sigma_1(x,t)+ \frac{1}{n_{\infty,2}}(1+ \frac{1-\varepsilon(1-\alpha)}{2} (d- m_2 |v|^2)) \sigma_2(x,t)   
\\&\quad + \frac{1}{n_{\infty,1}}\varepsilon (1- \delta) m_1 v \cdot \mu_1(x,t) + \frac{1}{n_{\infty,2}}(1-\frac{m_1}{m_2} \varepsilon (1- \delta)) m_2 v \cdot \mu_2(x,t) + \frac{1}{2} \frac{1}{n_{\infty,1}}\varepsilon (1- \alpha) (\frac{1}{d} m_2 |v|^2 -1) \tau_1(x,t) \\ &\quad +  \frac{1}{2} \frac{1}{n_{\infty,2}}(1-\varepsilon (1- \alpha)) (\frac{1}{d} m_2 |v|^2 -1) \tau_2(x,t)] -h_2\bigg). 
\end{split}
\end{align}}

In addition, we assume that $\sigma_k, \mu_k$ and $\tau_k$ are small enough such that we can neglect terms of the form $\sigma_k \sigma_l$, $\sigma_k \mu_l$ and $\sigma_k \tau_l$  ($k, l =1,2$). Thus we get 
{\footnotesize
\begin{align} \begin{split} \label{BGKhsmall}
&\quad\partial_t h_1 + v \cdot \nabla_x  h_1 \\ &= \nu_{11} n_{\infty,1}\left(f_1^{\infty} (v) [(\frac{1+d/2}{n_{\infty,1}} - \frac{m_1 |v|^2}{2 n_{\infty,1}}) \sigma_1(x,t) + \frac{m_1}{n_{\infty,1}} v \cdot \mu_1(x,t) + \frac{1}{n_{\infty,1}} (- \frac{1}{2} + \frac{m_1|v|^2}{2d}) \tau_1(x,t)]- h_1\right) \\&\quad + \nu_{12} n_{\infty,2}  (f_1^{\infty} [ \frac{1}{n_{\infty,1}}(1+ \frac{\alpha}{2} (d- m_1 |v|^2)) \sigma_1 + \frac{1}{2} \frac{1}{n_{\infty,2}}(1-\alpha) (d- m_1 |v|^2) \sigma_2 \\&\quad + \frac{1}{n_{\infty,1}}\delta m_1 v \cdot \mu_1 + \frac{1}{n_{\infty,2}}(1-\delta) m_1 v \cdot \mu_2 + \frac{1}{2} \frac{1}{n_{\infty,1}}\alpha (\frac{1}{d} m_1 |v|^2 -1) \tau_1 +  \frac{1}{2} \frac{1}{n_{\infty,2}}(1-\alpha) (\frac{1}{d} m_1 |v|^2 -1) \tau_2] - h_1),
\\[10pt]
&\quad\partial_t h_2 + v \cdot \nabla_x h_2 \\&  =\nu_{22} n_{\infty,2}\left(f_2^{\infty} (v) [(\frac{1+d/2}{n_{\infty,2}} - \frac{m_2 |v|^2}{2 n_{\infty,2}}) \sigma_2(x,t) + \frac{m_2}{n_{\infty,2}} v \cdot \mu_2(x,t) + \frac{1}{n_{\infty,2}} (- \frac{1}{2} + \frac{m_2|v|^2}{2d}) \tau_2(x,t)] -h_2\right)\\&\quad  + \nu_{21} n_{\infty,1} (f_2^{\infty} [\frac{1}{2} \frac{1}{n_{\infty,1}}\varepsilon(1-\alpha) (d- m_2 |v|^2) \sigma_1+ \frac{1}{n_{\infty,2}}(1+ \frac{1-\varepsilon(1-\alpha)}{2} (d- m_2 |v|^2)) \sigma_2  
 \\&\quad  + \frac{1}{n_{\infty,1}}\varepsilon (1- \delta) m_1 v \cdot \mu_1 + \frac{1}{n_{\infty,2}}(1-\frac{m_1}{m_2} \varepsilon (1- \delta)) m_2 v \cdot \mu_2 + \frac{1}{2} \frac{1}{n_{\infty,1}}\varepsilon (1- \alpha) (\frac{1}{d} m_2 |v|^2 -1) \tau_1 \\ &\quad +  \frac{1}{2} \frac{1}{n_{\infty,2}}(1-\varepsilon (1- \alpha)) (\frac{1}{d} m_2 |v|^2 -1) \tau_2] -h_2). 
\end{split}
\end{align}}

This is the linearized version of the system \eqref{BGK}. For this system, we want to prove exponential convergence to the equilibrium distributions \eqref{equ}. Explicitly, we want to prove the following results:
\begin{thm}
For each side length $L>0$ and dimension $d=1$, there exists an entropy functional $e(f_1, f_2)$ and a decay rate \textcolor{black}{$\tilde{C}$} satisfying 
$$ c_d(L)\, e(f_1, f_2) \leq ||f_1 - f_1^{\infty}||^2_{L^2\left((\frac{f_1^{\infty}(v)}{n_{\infty,1}})^{-1}dv d\tilde{x}\right)}+ ||f_2 - f_2^{\infty}||^2_{L^2\left((\frac{f_2^{\infty}(v)}{n_{\infty,2}})^{-1}dv d\tilde{x}\right)} \leq C_d(L)\, e(f_1, f_2)$$
with some positive constants $c_d$, $C_d$ \textcolor{black}{that depend on $L$.}

Moreover, assume that $$\nu_{11} n_{\infty,1} + \nu_{12} n_{\infty,2}= 1 \qquad\text{and }\qquad \nu_{22} n_{\infty,2} + \nu_{21} n_{\infty,1}= 1, $$ then any solution $(h_1(t), h_2(t))$ to \eqref{BGKhsmall} in 1D with $e(h_1(0) + f_1^{\infty}, h_2(0) + f_2^{\infty})< \infty,$ normalized according to \eqref{cons}, then satisfies
{\small $$e(h_1(t) + f_1^{\infty}, h_2(t)+ f_2^{\infty}) \leq e^{-\tilde{C} t}
e(h_1(0) + f_1^{\infty}, h_2(0) + f_2^{\infty}), $$}
where $\tilde{C}$ is given by {\scriptsize $$\tilde{C}=2 \min \lbrace\nu_{12} n_{\infty,2} (1- \delta), \nu_{12} n_{\infty,2} (1- \alpha), \nu_{11} n_{\infty,1}+ \nu_{12} n_{\infty,2}, \nu_{12} n_{\infty,1} \frac{m_1}{m_2} (1- \delta),  \nu_{12} n_{\infty,1} (1- \alpha), \nu_{22} n_{\infty,2} + \nu_{12} n_{\infty,1} , 2 \mu\rbrace.$$ }
Here, \textcolor{black}{$\mu^d(L)$} is a one species decay rate developped in \textcolor{black}{theorem 1.1 in } \cite{achleitner2017multi}. 

\label{theo}
\end{thm}
\subsection{Linearized BGK equation in 1D}
In this section we consider the system of linearized BGK equations for gas mixtures in 1D. For this system, we want to analyse the large time behaviour in the following section. The idea is to use stability criteria for ODEs. Therefore, we want to transform the system \eqref{BGKhsmall} to a system of infinite ODEs. 
Compared to the one species case, we expand each equation in \eqref{BGKhsmall} in a different Hilbert space, namely 
a different weighted $L^2$ space in the velocity. The system of equations in \eqref{BGKhsmall} in 1D is given by
\begin{align} \begin{split} \label{BGKh1D}
\scriptstyle
\partial_t h_1 + v \cdot \nabla_x  h_1  &\scriptstyle= \nu_{11} n_{\infty,1}\left(f_1^{\infty} (v) [(\frac{1+1/2}{n_{\infty,1}} - \frac{m_1 v^2}{2 n_{\infty,1}}) \sigma_1(x,t) + \frac{m_1}{n_{\infty,1}} v \cdot \mu_1(x,t) + \frac{1}{n_{\infty,1}} (- \frac{1}{2} + \frac{m_1 v^2}{2}) \tau_1(x,t)]- h_1\right)\\&\scriptstyle\quad+ \nu_{12} n_{\infty,2} \big(f_1^{\infty} [ \frac{1}{n_{\infty,1}}(1+ \frac{\alpha}{2} (1- m_1 v^2)) \sigma_1 + \frac{1}{2} \frac{1}{n_{\infty,2}}(1-\alpha) (1- m_1 v^2) \sigma_2 
\\&\scriptstyle\quad+ \frac{1}{n_{\infty,1}}\delta m_1 v \cdot \mu_1 + \frac{1}{n_{\infty,2}}(1-\delta) m_1 v \cdot \mu_2 + \frac{1}{2} \frac{1}{n_{\infty,1}}\alpha ( m_1 v^2 -1) \tau_1 +  \frac{1}{2} \frac{1}{n_{\infty,2}}(1-\alpha) ( m_1 v^2 -1) \tau_2] - h_1\big),
\\ 
\scriptstyle\partial_t h_2 + v \cdot \nabla_x h_2 &=\scriptstyle\nu_{22} n_{\infty,2}\left(f_2^{\infty} (v) [(\frac{1+1/2}{n_{\infty,2}} - \frac{m_2 v^2}{2 n_{\infty,2}}) \sigma_2(x,t) + \frac{m_2}{n_{\infty,2}} v \cdot \mu_2(x,t) + \frac{1}{n_{\infty,2}} (- \frac{1}{2} + \frac{m_2 v^2}{2}) \tau_2(x,t)] -h_2\right)
\\& \scriptstyle\quad+ \nu_{21} n_{\infty,1} \big(f_2^{\infty} [\frac{1}{2} \frac{1}{n_{\infty,1}}\varepsilon(1-\alpha) (1- m_2 v^2) \sigma_1+ \frac{1}{n_{\infty,2}}(1+ \frac{1-\varepsilon(1-\alpha)}{2} (1- m_2 v^2)) \sigma_2   \\&\scriptstyle\quad+ \frac{1}{n_{\infty,1}}\varepsilon (1- \delta) m_1 v \cdot \mu_1 + \frac{1}{n_{\infty,2}}(1-\frac{m_1}{m_2} \varepsilon (1- \delta)) m_2 v \cdot \mu_2 + \frac{1}{2} \frac{1}{n_{\infty,1}}\varepsilon (1- \alpha) ( m_2 v^2 -1) \tau_1
 \\ &\scriptstyle\quad +\frac{1}{2} \frac{1}{n_{\infty,2}}(1-\varepsilon (1- \alpha)) ( m_2 v^2 -1) \tau_2] -h_2 \big),  
\end{split}
\end{align}
for the perturbations $$h_1(x,v,t) \approx f_1(x,v,t) - f_1^{\infty}, \qquad h_2(x,v,t) \approx f_2(x,v,t) - f_2^{\infty}. $$

In order to get rid of the derivatives in $x$-space, we expand $h_1,h_2$ in the $x$-Fourier series
$$ h_1(x,v,t) = \sum_{k \in \mathbb{Z}} h_{1,k}(v,t) e^{ik \frac{2 \pi}{L} x}, \quad h_2(x,v,t) = \sum_{k \in \mathbb{Z}} h_{2,k}(v,t) e^{ik \frac{2 \pi}{L} x}. $$
We insert this expansion into \eqref{BGKh1D} and obtain for each spatial mode $h_{1,k}(v,t)$ and $h_{2,k}(v,t)$
\begin{align} \begin{split} \label{BGKh1Dmodes} &\scriptstyle
\quad\partial_t h_{1,k} + i k \frac{2 \pi}{L} v  h_{1,k} \\  &\scriptstyle= \nu_{11} n_{\infty,1}\left(f_1^{\infty} (v) [(\frac{1+1/2}{n_{\infty,1}} - \frac{m_1 v^2}{2 n_{\infty,1}}) \sigma_{1,k}(t) + \frac{m_1}{n_{\infty,1}} v \cdot \mu_{1,k}(t) + \frac{1}{n_{\infty,1}} (- \frac{1}{2} + \frac{m_1 v^2}{2}) \tau_{1,k}(t)]- h_{1,k}\right)
\\&\scriptstyle\quad+ \nu_{12} n_{\infty,2} \bigg(f_1^{\infty}(v) [ \frac{1}{n_{\infty,1}}(1+ \frac{\alpha}{2} (1- m_1 v^2)) \sigma_{1,k}(t) + \frac{1}{2} \frac{1}{n_{\infty,2}}(1-\alpha) (1- m_1 v^2) \sigma_{2,k}(t) 
\\&\scriptstyle\quad+ \frac{1}{n_{\infty,1}}\delta m_1 v \cdot \mu_1(t) + \frac{1}{n_{\infty,2}}(1-\delta) m_1 v \cdot \mu_{2,k}(t) + \frac{1}{2} \frac{1}{n_{\infty,1}}\alpha ( m_1 v^2 -1) \tau_1(t) +  \frac{1}{2} \frac{1}{n_{\infty,2}}(1-\alpha) ( m_1 v^2 -1) \tau_{2,k}(t)] - h_{1,k}\bigg), 
\\[8pt]&\scriptstyle\quad\partial_t h_{2,k} + i k \frac{2 \pi}{L} v h_{2,k} \\&\scriptstyle=\nu_{22} n_{\infty,2}\left(f_2^{\infty} (v) [(\frac{1+1/2}{n_{\infty,2}} - \frac{m_2 v^2}{2 n_{\infty,2}}) \sigma_{2,k}(t) + \frac{m_2}{n_{\infty,2}} v \cdot \mu_{2,k}(t) + \frac{1}{n_{\infty,2}} (- \frac{1}{2} + \frac{m_2 v^2}{2}) \tau_{2,k}(t)] -h_{2,k}\right)\\&\scriptstyle\quad + \nu_{21} n_{\infty,1}\bigg(f_2^{\infty} [\frac{1}{2} \frac{1}{n_{\infty,1}}\varepsilon(1-\alpha) (1- m_2 v^2) \sigma_{1,k}(t)+ \frac{1}{n_{\infty,2}}(1+ \frac{1-\varepsilon(1-\alpha)}{2} (1- m_2 v^2)) \sigma_{2,k}(t)   \\&\scriptstyle\quad + \frac{1}{n_{\infty,1}}\varepsilon (1- \delta) m_1 v \cdot \mu_{1,k}(t) + \frac{1}{n_{\infty,2}}(1-\frac{m_1}{m_2} \varepsilon (1- \delta)) m_2 v \cdot \mu_{2,k}(t) + \frac{1}{2} \frac{1}{n_{\infty,1}}\varepsilon (1- \alpha) ( m_2 v^2 -1) \tau_{1,k}(t) \\ &\scriptstyle\quad +\frac{1}{2} \frac{1}{n_{\infty,2}}(1-\varepsilon (1- \alpha)) ( m_2 v^2 -1) \tau_{2,k}(t)] -h_{2,k}\bigg),  
\end{split}
\end{align}
where $\sigma_{1,k}, \sigma_{2,k}, \mu_{1,k}, \mu_{2,k}, \tau_{1,k}$ and $\tau_{2,k}$ denote the spatial modes of the moments $\sigma_{1}, \sigma_{2}, \mu_{1}, \mu_{2}, \tau_{1}$ and $\tau_{2}$ given by
$$ \sigma_{1,k} = \int h_{1,k}(v,t) dv,~ \sigma_{2,k} = \int h_{2,k}(v,t) dv,$$ $$ \mu_{1,k} = \int v h_{1,k}(v,t) dv, ~ \mu_{2,k} = \int v h_{2,k}(v,t) dv,$$ $$ \tau_{1,k} = \int m_1 |v|^2 h_{1,k}(v,t) dv, ~ \tau_{2,k} = \int m_2 |v|^2 h_{2,k}(v,t) dv.$$
Now, define the functions $g_{1,0}, g_{1,1}, g_{1,2}, g_{2,0}, g_{2,1}$ and $g_{2,2}$ as
\begin{align*}
g_{1,0}(v) = f_1^{\infty}(v) \frac{1}{n_{\infty,1}},\quad  g_{1,1}(v) =\sqrt{m_1} v f_1^{\infty}(v) \frac{1}{n_{\infty,1}} ,\quad g_{1,2}(v) = \frac{ m_1 v^2 -1}{\sqrt{2}}  f_1^{\infty}(v) \frac{1}{n_{\infty,1}}, \\ g_{2,0}(v) = f_2^{\infty}(v) \frac{1}{n_{\infty,2}}, \quad g_{2,1}(v) = \sqrt{m_2} v f_2^{\infty}(v) \frac{1}{n_{\infty,2}}, \quad g_{2,2}(v) = \frac{m_2 v^2 - 1}{\sqrt{2}} f_2^{\infty}(v) \frac{1}{n_{\infty,2}}. 
\end{align*}
Then we can rewrite equation \eqref{BGKh1D} as
\begin{align} \begin{split} \label{BGKh1Dnew}
\scriptstyle\partial_t h_{1,k} + i k \frac{2 \pi}{L} v  h_{1,k}  &\scriptstyle= \nu_{11} n_{\infty,1}  ( g_{1,0}(v) \sigma_{1,k}(t) + \sqrt{m_1} g_{1,1}(v) \mu_{1,k}(t) + g_{1,2}(v) \frac{1}{\sqrt{2}}(\tau_{1,k}(t) - \sigma_{1,k}(t))- h_{1,k}) \\&\scriptstyle\quad+ \nu_{12} n_{\infty,2}  (g_{1,0}(v) \sigma_{1,k}(t) + \delta \sqrt{m_1} g_{1,1}(v) \mu_{1,k}(t) + (1- \delta) \frac{n_{\infty,1}}{n_{\infty,2}} \sqrt{m_1} g_{1,1}(v) \mu_{2,k}(t) \\&\scriptstyle\quad+ \alpha g_{1,2}(v) \frac{1}{\sqrt{2}} (\tau_{1,k}(t) - \sigma_{1,k}(t)) + (1- \alpha) \frac{n_{\infty,1}}{n_{\infty,2}} g_{1,2}(v) \frac{1}{\sqrt{2}} (\tau_{2,k}(t) - \sigma_{2,k}(t)) - h_{1,k}),
\\[8pt]
\scriptstyle\partial_t h_{2,k} + i k \frac{2 \pi}{L} v h_{2,k} &\scriptstyle=\nu_{22} n_{\infty,2}  (g_{2,0}(v) \sigma_{2,k}(t) + \sqrt{m_2} g_{2,1}(v) \mu_{2,k}(t) + g_{2,2}(v) \frac{1}{\sqrt{2}} (\tau_{2,k}(t) - \sigma_{2,k}(t)) -h_{2,k})\\& \scriptstyle\quad+ \nu_{21} n_{\infty,1}  (g_{2,0}(v) \sigma_{2,k}(t) + \frac{n_{\infty,2}}{n_{\infty,1}} \frac{m_1}{m_2} \varepsilon (1- \delta) \sqrt{m_2} g_{2,1}(v) \mu_{1,k}(t) + (1- \frac{m_1}{m_2} \varepsilon (1- \delta)) \sqrt{m_2} g_{2,1}(v) \mu_{2,k}(t) \\&\scriptstyle\quad+\frac{n_{\infty,2}}{n_{\infty,1}} \varepsilon (1- \alpha) g_{2,2}(v) \frac{1}{\sqrt{2}} (\tau_{1,k}(t) - \sigma_{1,k}(t)) + (1- \varepsilon (1- \alpha)) \frac{1}{\sqrt{2}} g_{2,2}(v) (\tau_{2,k}(t) - \sigma_{2,k}(t)) -h_{2,k}). 
\end{split}
\end{align}
Note that the functions $g_{1,0}, g_{1,1}, g_{1,2}$ satisfy 
\begin{align*}
&\int g_{1,0}(v) g_{1,0}(v) \left( \frac{f_1^{\infty}(v)}{n_{\infty,1}} \right)^{-1} dv = \int f_1^{\infty}(v) \frac{1}{n_{\infty,1}} dv = 1 \\
&\int g_{1,1}(v) g_{1,1}(v) \left( \frac{f_1^{\infty}(v)}{n_{\infty,1}} \right)^{-1} dv = \int f_1^{\infty}(v) \frac{1}{n_{\infty,1}} m_1 v^2 dv = 1 \\
&\int g_{1,2}(v) g_{1,2}(v) \left( \frac{f_1^{\infty}(v)}{n_{\infty,1}} \right)^{-1} dv = \int \frac{m_1 v^2 -1}{\sqrt{2}} f_1^{\infty}(v) \frac{1}{n_{\infty,1}} \frac{m_1 v^2 -1}{\sqrt{2}} dv = 1 \\
&\int g_{1,0}(v) g_{1,1}(v) \left( \frac{f_1^{\infty}(v)}{n_{\infty,1}} \right)^{-1} dv = \int v f_1^{\infty}(v) \frac{1}{n_{\infty,1}} \sqrt{m_1} dv = 0 \\
&\int g_{1,0}(v) g_{1,2}(v) \left( \frac{f_1^{\infty}(v)}{n_{\infty,1}} \right)^{-1} dv = \int \frac{m_1 v^2 -1}{\sqrt{2}} f_1^{\infty}(v) \frac{1}{n_{\infty,1}}  dv = 0 \\
&\int g_{1,1}(v) g_{1,2}(v) \left( \frac{f_1^{\infty}(v)}{n_{\infty,1}} \right)^{-1} dv = \int v \sqrt{m_1} \frac{m_1 v^2 -1}{\sqrt{2}} f_1^{\infty}(v) \frac{1}{n_{\infty,1}}dv = 0. 
\end{align*}
In the same way one can prove that $g_{2,0}, g_{2,1}$ and $g_{2,2}$ are orthonormal in $L^2 \left(\mathbb{R}; (\frac{f_2^{\infty}(v)}{n_{\infty,2}})^{-1}\right)$.

Now we extend $g_{1,0}, g_{1,1}, g_{1,2}$ to an orthonormal basis $\lbrace g_{1,m} (v) \rbrace_{m \in \mathbb{N}_0}$ in $L^2 \left(\mathbb{R}; (\frac{f_1^{\infty}(v)}{n_{\infty,1}})^{-1}\right)$ and $g_{2,0},$ $ g_{2,1},$ $ g_{2,2}$ to an orthonormal basis $\lbrace g_{2,m} (v) \rbrace_{m \in \mathbb{N}_0}$ in $L^2 \left(\mathbb{R}; (\frac{f_2^{\infty}(v)}{n_{\infty,2}})^{-1}\right)$. 
One can expand $h_{1,k}( \cdot , t) \in L^2 \left(\mathbb{R}; (\frac{f_1^{\infty}(v)}{n_{\infty,1}})^{-1}\right)$ and $h_{2,k}( \cdot , t) \in L^2 \left(\mathbb{R};(\frac{f_2^{\infty}(v)}{n_{\infty,2}})^{-1}\right)$ in the corresponding orthonormal basis
\begin{align}
\begin{split}
h_{1,k}(v,t) = \sum_{m=0}^{\infty} \hat{h}_{1,(k,m)} g_{1,m}(v) \quad \text{with} \quad \hat{h}_{1,(k,m)} = \langle h_{1,k}(v), g_{1,m}(v) \rangle_{L^2 \left((\frac{f_1^{\infty}(v)}{n_{\infty,1}})^{-1}\right)} \\
h_{2,k}(v,t) = \sum_{m=0}^{\infty} \hat{h}_{2,(k,m)} g_{2,m}(v) \quad \text{with} \quad \hat{h}_{2,(k,m)} = \langle h_{2,k}(v), g_{2,m}(v) \rangle_{L^2 \left((\frac{f_2^{\infty}(v)}{n_{\infty,2}})^{-1}\right)}. 
\end{split}
\label{expg}
\end{align}
For each $k \in \mathbb{Z}$, the infinite vectors
\begin{align}
\label{InfVec}
\begin{split}
&\hat{h}_{1,k}(t) = (\hat{h}_{1,(k,0)}(t), \,\hat{h}_{1,(k,1)}(t), \cdots )^T \in l^2(\mathbb{N}_0), \\[2pt]
&\hat{h}_{2,k}(t) = (\hat{h}_{2,(k,0)}(t), \,\hat{h}_{2,(k,1)}(t), \cdots )^T \in l^2(\mathbb{N}_0)
\end{split}
\end{align}
contain all the coefficients of $h_{1,k}(\cdot,t)$ and $h_{2,k}(\cdot,t)$ in the expansion \eqref{expg}, respectively. 

In particular, one has
\begin{align*}
&\hat{h}_{1,(k,0)} = \int h_{1,k}(v) g_{1,0}(v) \left( \frac{f_1^{\infty}(v)}{n_{\infty,1}} \right)^{-1} dv = \int h_{1,k}(v) dv = \sigma_{1,k}, \\
&\hat{h}_{1,(k,1)} = \int h_{1,k}(v) g_{1,1}(v) \left( \frac{f_1^{\infty}(v)}{n_{\infty,1}} \right)^{-1} dv = \int h_{1,k}(v) \sqrt{m_1} v dv = \sqrt{m_1} \mu_{1,k}, \\
&\hat{h}_{1,(k,2)} = \int h_{1,k}(v) g_{1,2}(v) \left( \frac{f_1^{\infty}(v)}{n_{\infty,1}} \right)^{-1} dv = \int h_{1,k}(v) \frac{m_1 v^2 - 1}{\sqrt{2}} dv = \frac{1}{\sqrt{2}} ( \tau_{1,k} - \sigma_{1,k}), \\
&\hat{h}_{2,(k,0)} = \int h_{2,k}(v) g_{2,0}(v) \left( \frac{f_2^{\infty}(v)}{n_{\infty,2}} \right)^{-1} dv = \int h_{2,k}(v) dv = \sigma_{2,k},  \\
&\hat{h}_{2,(k,1)} = \int h_{2,k}(v) g_{2,1}(v) \left( \frac{f_2^{\infty}(v)}{n_{\infty,2}} \right)^{-1} dv = \int h_{2,k}(v) \sqrt{m_2} v dv = \sqrt{m_2} \mu_{2,k},  \\
&\hat{h}_{2,(k,2)} = \int h_{2,k}(v) g_{2,2}(v) \left( \frac{f_2^{\infty}(v)}{n_{\infty,2}} \right)^{-1} dv = \int h_{2,k}(v) \frac{m_2 v^2 - 1}{\sqrt{2}} dv = \frac{1}{\sqrt{2}} (\tau_{2,k} - \sigma_{2,k}). 
\end{align*}
Hence, \eqref{BGKh1Dnew} can be written equivalently as 
{\small
\begin{align} \begin{split} \label{BGKh1Dnew2}
\partial_t h_{1,k} + i k \frac{2 \pi}{L} v  h_{1,k}  &= \nu_{11} n_{\infty,1}  ( g_{1,0}(v) \hat{h}_{1,(k,0)}(t) +  g_{1,1}(v) \hat{h}_{1,(k,1)}(t) + g_{1,2}(v) \hat{h}_{1,(k,2)}(t)- h_{1,k}) \\&\quad + \nu_{12} n_{\infty,2}  (g_{1,0}(v) \hat{h}_{1,(k,0)} + \delta  g_{1,1}(v) \hat{h}_{1,(k,1)} + (1- \delta) \frac{n_{\infty,1}}{n_{\infty,2}} \frac{\sqrt{m_1}}{\sqrt{m_2}} g_{1,1}(v) \hat{h}_{2,(k,1)}(t) \\&\quad +\alpha g_{1,2}(v) \hat{h}_{1,(k,2)}(t)+(1- \alpha) \frac{n_{\infty,1}}{n_{\infty,2}} g_{1,2}(v) \hat{h}_{2,(k,2)}(t) - h_{1,k}),
\\[6pt]
\partial_t h_{2,k} + i k \frac{2 \pi}{L} v h_{2,k} &=\nu_{22} n_{\infty,2}  (g_{2,0}(v) \hat{h}_{2,(k,0)}(t) +  g_{2,1}(v) \hat{h}_{2,(k,1)}(t) + g_{2,2}(v) \hat{h}_{2,(k,2)}(t) -h_{2,k})\\&\quad + \nu_{21} n_{\infty,1}  (g_{2,0}(v) \hat{h}_{2,(k,0)}(t) + \frac{n_{\infty,2}}{n_{\infty,1}} \frac{\sqrt{m_1}}{\sqrt{m_2}} \varepsilon (1- \delta)  g_{2,1}(v) \hat{h}_{1,(k,1)}(t) \\&\quad + (1- \frac{m_1}{m_2} \varepsilon (1- \delta))  g_{2,1}(v) \hat{h}_{2,(k,1)}(t) + \frac{n_{\infty,2}}{n_{\infty,1}} \varepsilon (1- \alpha) g_{2,2}(v) \hat{h}_{1,(k,2)}(t)\\&\quad + (1- \varepsilon (1- \alpha)) \frac{1}{\sqrt{2}} g_{2,2}(v) \hat{h}_{2,(k,2)}(t) -h_{2,k}). 
\end{split}
\end{align}}
Therefore, by using \eqref{expg} and conducting projection onto the corresponding weighted $L^2$ in velocity space, 
one gets that the vectors of coefficients defined in (\ref{InfVec}) satisfy
\begin{align}
\begin{split}
\frac{d}{dt} \hat{h}_{1,k}(t) + i k \frac{2\pi}{L} L_{1,1} \hat{h}_{1,k}(t) = - \nu_{11} n_{\infty,1} L_{1,2} \hat{h}_{1,k}(t) - \nu_{12} n_{\infty,2} L_{1,3} \hat{h}_{1,k}(t) + \nu_{12} n_{\infty,2} L_{1,4} \hat{h}_{2,k}(t), \\[4pt]
\frac{d}{dt} \hat{h}_{2,k}(t) + i k \frac{2\pi}{L} L_{2,1} \hat{h}_{2,k}(t) = - \nu_{22} n_{\infty,2} L_{2,2} \hat{h}_{2,k}(t) - \nu_{21} n_{\infty,1} L_{2,3} \hat{h}_{2,k}(t) + \nu_{21} n_{\infty,1} L_{2,4} \hat{h}_{1,k}(t), 
\end{split}
\label{BGKend}
\end{align}
where $L_{1,1}, L_{1,2}, L_{1,3}, L_{1,4}, L_{2,1}, L_{2,2}, L_{2,3}$ and $L_{2,4}$ are represented by ``infinite matrices" on \\ $l^2(\mathbb{N}_0)$ given by
\begin{align*}
 &L_{1,1}=L_{2,1}=\begin{pmatrix}
0 & \sqrt{1} & 0 & \cdots \\ \sqrt{1} & 0 & \sqrt{2} & 0 \\ 0 & \sqrt{2} & 0 & \sqrt{3} \\ \vdots & 0 & \sqrt{3} & \ddots
\end{pmatrix},
\quad
L_{1,2}= L_{2,2} = \text{diag}(0,0,0,1,1, \cdots),  \\[4pt]
&L_{1,3} = \text{diag}(0, (1- \delta), (1- \alpha), 1,1, \cdots ), \\[4pt]
&L_{1,4} = \text{diag}(0, (1- \delta) \frac{n_{\infty,1}}{n_{\infty,2}} \frac{\sqrt{m_1}}{\sqrt{m_2}}, (1- \alpha) \frac{n_{\infty,1}}{n_{\infty,2}}, 0, 0, \cdots ), \\[4pt]
&L_{2,3} = \text{diag}(0, \frac{m_1}{m_2} \varepsilon (1- \delta), \varepsilon(1- \alpha), 1,1, \cdots ), \\[4pt]
&L_{2,4} = \text{diag}(0, \frac{n_{\infty,2}}{n_{\infty,1}} \frac{\sqrt{m_1}}{\sqrt{m_2}} \varepsilon (1- \delta) , \frac{n_{\infty,2}}{n_{\infty,1}} \varepsilon (1- \alpha) , 0, 0, \cdots). 
\end{align*}
Note that $L_{1,1}, L_{1,2}, L_{1,3}, L_{1,4}$ represent coefficients with respect to a different basis than $L_{2,1},$ $ L_{2,2,},$ $ L_{2,3},$ $ L_{2,4}$. As a consequence, for example, $L_{1,4} \hat{h}_{2,k}(t)$ has a different meaning than $L_{2,3} \hat{h}_{2,k}(t)$ even if $L_{1,4}= L_{2,3}$. 

In a word, we obtained a system of infinite ordinary differential equations  given by (\ref{BGKend}). 

\section{Hypocoercivity estimate}
In this section we want to prove the estimate stated in theorem \ref{theo}.
\subsection{Definition of the entropy functional}
For the definition of the entropy functional for the gas mixture, we take the natural choice from a physical point of view and simply take a weighted sum of 
the entropies of species 1 and species 2. 

We consider a solution $(h_1,h_2)$ of \eqref{BGKh2} and define the entropy functional for the gas mixture  entropy functional $e(\tilde{f}_1, \tilde{f}_2)$ by
\begin{align}
e(\tilde{f}_1,\tilde{f}_2):=  \sum_{k \in \mathbb{Z}} \left( \frac{1}{n_{\infty,1}}\langle h_{1,k}(v), P_k h_{1,k}(v) \rangle_{L^2(( \frac{f_1^{\infty}(v)}{n_{\infty,1}})^{-1})} + \frac{1}{n_{\infty,2}} \langle h_{2,k}(v), P_k h_{2,k}(v)\rangle_{L^2(( \frac{f_2^{\infty}(v)}{n_{\infty,2}})^{-1})} \right) 
\label{entropy}
\end{align}
with $$\tilde{f}_1 (t) = f_1^{\infty} + h_1(t), \qquad\tilde{f}_2(t) = f_2^{\infty} + h_2(t)$$ and the 
``infinite matrices" $P_0=\textbf{1}$ and $P_k, ~k>0$ from the one species case having
$$\begin{pmatrix}
1 & - \frac{i \tilde{\alpha}}{k} & 0 & 0 \\ i \frac{\tilde{\alpha}}{k} & 1 & - \frac{i \beta}{k} & 0 \\ 0 & \frac{i \beta }{k} & 1 & - \frac{i \tilde{\gamma}}{k} \\ 0 & 0 & \frac{i \tilde{\gamma}}{k} & 1 
\end{pmatrix}$$
with $\tilde{\alpha}>0, $ $ \beta>0$ and $\tilde{\gamma} >0$ specified later, as upper left $4\times 4$ block with all other entries being those of the identity. For details of determining this matrix in the one species case see \cite{achleitner2017multi}. Here, the infinite matrices $P_0$ and $P_k$ for $k>0$ are regarded as bounded operators in $L^2(( \frac{f_1^{\infty}(v)}{n_{\infty,1}})^{-1})$ in the first term in the entropy and in $L^2(( \frac{f_2^{\infty}(v)}{n_{\infty,2}})^{-1})$ in the second term in the entropy. 
\textcolor{black}{
\begin{rmk}
Note that the entropy defined in \eqref{entropy} is not the only possible choice. For example, the analysis also works through if defining the entropy as
\begin{align*}
e(\tilde{f}_1,\tilde{f}_2):=  \sum_{k \in \mathbb{Z}} \left( \frac{n_{\infty,2}}{n_{\infty,1}+n_{\infty,2}}\langle h_{1,k}(v), P_k h_{1,k}(v) \rangle_{L^2(( \frac{f_1^{\infty}(v)}{n_{\infty,1}})^{-1})} + \frac{n_{\infty,1}}{n_{\infty,1}+n_{\infty,2}} \langle h_{2,k}(v), P_k h_{2,k}(v)\rangle_{L^2(( \frac{f_2^{\infty}(v)}{n_{\infty,2}})^{-1})} \right). 
\end{align*}
\label{rem}
\end{rmk}}

We now insert the expansions \eqref{expg} in this total entropy and obtain
\begin{align*}
e(\tilde{f}_1,\tilde{f}_2)= \sum_{k \in \mathbb{Z}}\Bigg[\frac{1}{n_{\infty,1}} \langle \sum_{m=0}^{\infty} \hat{h}_{1,(k,m)}(t) g_{1,m}(v), \sum_{l=0}^{\infty} (P_k \hat{h}_{1,k})_{(l)} g_{1,l}(v)\rangle_{L^2(( \frac{f_1^{\infty}(v)}{n_{\infty,1}})^{-1})}\\[2pt]\qquad + \frac{1}{n_{\infty,2}} \langle \sum_{m=0}^{\infty} \hat{h}_{2,(k,m)}(t) g_{2,m}(v), \sum_{l=0}^{\infty} (P_k \hat{h}_{2,k})_{(l)} g_{2,l}(v)\rangle_{L^2(( \frac{f_2^{\infty}(v)}{n_{\infty,2}})^{-1})}\Bigg], 
\end{align*}
where $(P_k \hat{h}_{1,k})_{(l)}$ denotes the $l$-th component of $P_k \hat{h}_{1,k}$ in the expansion in $L^2(( \frac{f_1^{\infty}(v)}{n_{\infty,1}})^{-1})$ with respect to $\lbrace g_{1,m}(v)\rbrace_{m \in \mathbb{N}_0}$ and $(P_k \hat{h}_{2,k})_{(l)}$ denotes the $l$-th component of $P_k \hat{h}_{2,k}$ in the expansion in $L^2(( \frac{f_2^{\infty}(v)}{n_{\infty,2}})^{-1})$ with respect to $\lbrace g_{2,m}(v)\rbrace_{m \in \mathbb{N}_0}$. 

If we compute the Cauchy product of the two rows, we get
\begin{align*}
e(\tilde{f}_1,\tilde{f}_2)= \sum_{k \in \mathbb{Z}} \Bigg[\frac{1}{n_{\infty,1}}\sum_{m=0}^{\infty} \sum_{l=0}^m \overline{\hat{h}_{1,(k,m-l)}} (P_k \hat{h}_{1,k})_{(l)} \langle g_{1,(m-l)}(v), g_{1,l}(v) \rangle_{L^2(( \frac{f_1^{\infty}(v)}{n_{\infty,1}})^{-1})} \\[2pt]\qquad + \frac{1}{n_{\infty,2}}\sum_{m=0}^{\infty} \sum_{l=0}^m \overline{\hat{h}_{2,(k,m-l)}} (P_k \hat{h}_{2,k})_{(l)} \langle g_{2,(m-l)}(v), g_{2,l}(v) \rangle_{L^2(( \frac{f_2^{\infty}(v)}{n_{\infty,2}})^{-1})}\Bigg]. 
\end{align*}
The bar denotes the complex conjugate and comes from the scalar product in complex space. Since $\lbrace g_{1,m}(v)\rbrace_{m \in \mathbb{N}_0}$ is orthonormal in $L^2(( \frac{f_1^{\infty}(v)}{n_{\infty,1}})^{-1})$, and $\lbrace g_{2,m}(v)\rbrace_{m \in \mathbb{N}_0}$ is orthonormal in \\$L^2(( \frac{f_2^{\infty}(v)}{n_{\infty,2}})^{-1})$, we obtain
{\footnotesize
\begin{align*}
e(\tilde{f}_1,\tilde{f}_2)&= \sum_{k \in \mathbb{Z}}\left( \frac{1}{n_{\infty,1}} \sum_{m=0}^{\infty} \sum_{l=0}^m \overline{\hat{h}_{1,(k,m-l)}} (P_k \hat{h}_{1,k})_{(l)} \delta_{(m-l,l)} + \frac{1}{n_{\infty,2}} \sum_{m=0}^{\infty} \sum_{l=0}^m \overline{\hat{h}_{2,(k,m-l)}} (P_k \hat{h}_{2,k})_{(l)}  \delta_{(m-l,l)}\right). 
\end{align*}}
We can rewrite the term 
\begin{align}
e_{k,1}(\tilde{f}_1):=\sum_{m=0}^{\infty} \sum_{l=0}^m \overline{\hat{h}_{1,(k,m-l)}} (P_k \hat{h}_{1,k})_{(l)} \delta_{(m-l,l)}
\label{e_k}
\end{align}
 as 
\begin{align*}
\lim_{S \rightarrow \infty} \sum_{\substack{m=0 \\ \text{m even}}}^S \overline{\hat{h}_{1,(k, \frac{m}{2}})} (P_k \hat{h}_{1,k})_{(\frac{m}{2})} = \lim_{S \rightarrow \infty} \sum_{l=0}^{\frac{S}{2}} \overline{\hat{h}_{1,(k,l)}} (P_k \hat{h}_{1,k})_{(l)} = \lim_{M \rightarrow \infty} \sum_{l=0}^M \overline{\hat{h}_{1,(k,l)}} (P_k \hat{h}_{1,k})_{(l)}\\[2pt]
 = \lim_{M \rightarrow \infty} \hat{h}^{(M)}_{1,k} \cdot (P_k \hat{h}_{1,k})^{(M)} = \lim_{M \rightarrow \infty} \hat{h}_{1,k}^{(M)} \cdot P_k^{(M \times M)} \cdot \hat{h}_{1,k}^{(M)}, 
\end{align*}
 where the upper index $(M)$ indicates that we take an $(M+1)$-dimensional vector with the first $M+1$ entries of the corresponding ``infinite vector", 
 defined as the vector containing the coefficients of the expansion similar to \eqref{expg}. 
Similar for the upper index $(M \times M)$. Here we take the upper $(M+1) \times (M+1)$ left block of the corresponding ``infinite matrix". The dots are the notation for the scalar product, meaning $\hat{h} \cdot P \cdot \hat{h}=\overline{\hat{h}^{T}} P \hat{h}$ with a vector $\hat{h}$ and matrix $P$. 

In the same way, we get for the second species term
\begin{align*}
e_{k,2}(\tilde{f}_2)=  \sum_{m=0}^{\infty} \sum_{l=0}^m \overline{\hat{h}_{2,(k,m-l)}} (P_k \hat{h}_{2,k})_{(l)} \delta_{(m-l,l)}= \lim_{M \rightarrow \infty} \hat{h}_{2,k}^{(M)} \cdot P_k^{(M \times M)} \cdot \hat{h}_{2,k}^{(M)}.
 \end{align*}
 
Now, we want to consider the time derivative of the total entropy $e(f_1,f_2)$. We have that
\begin{align*}
 \frac{d}{dt} (\frac{1}{n_{\infty,1}} e_{k,1}(\tilde{f}_1) + \frac{1}{n_{\infty,2}} e_{k,2}(\tilde{f}_2))= \frac{1}{n_{\infty,1}} \lim_{M \rightarrow \infty} ( \frac{d}{dt}\hat{h}_{1,k}^{(M)} \cdot P_k^{(M \times M)} \cdot \hat{h}_{1,k}^{(M)} + \hat{h}_{1,k}^{(M)} \cdot P_k^{(M \times M)} \cdot \frac{d}{dt}\hat{h}_{1,k}^{(M)})\\[4pt]\qquad +\frac{1}{n_{\infty,2}} \lim_{M \rightarrow \infty} (\frac{d}{dt}\hat{h}_{2,k}^{(M)} \cdot P_k^{(M \times M)} \cdot \hat{h}_{2,k}^{(M)} + \hat{h}_{2,k}^{(M)} \cdot P_k^{(M \times M)} \cdot \frac{d}{dt}\hat{h}_{2,k}^{(M)}). 
 \end{align*}
We want to show that we can estimate the right-hand side by $e(\tilde{f}_1,\tilde{f}_2)$ in order to get an estimate for $e(\tilde{f}_1,\tilde{f}_2)$ using the Gronwall's estimate. 
 
For this, we need estimates on
{\small
$$\textcolor{black}{E_k:=}\frac{1}{n_{\infty,1}} \Big(\hat{h}_{1,k}^{(M)} \cdot P_k^{(M \times M)} \cdot \frac{d}{dt}\hat{h}_{1,k}^{(M)}+\frac{d}{dt}\hat{h}_{1,k}^{(M)} \cdot P_k^{(M \times M)} \cdot \hat{h}_{1,k}^{(M)}\Big)+\frac{1}{n_{\infty,2}}\Big(\hat{h}_{2,k}^{(M)} \cdot P_k^{(M \times M)} \cdot \frac{d}{dt}\hat{h}_{2,k}^{(M)}+
\frac{d}{dt}\hat{h}_{2,k}^{(M)} \cdot P_k^{(M \times M)} \cdot \hat{h}_{2,k}^{(M)}\Big),$$}
which will be derived in the following. In comparison to one species, we start with a system of two partial differential equations which are coupled due to interspecies interactions. This requires additional estimates on these interspecies terms. These estimates are very different from the one species case, because the interspecies terms from species 1 and the interspecies terms from species 2 
have to be coupled and estimated as a whole. Therefore, we need more delicate techniques.
 
\subsection{Estimates on \textcolor{black}{$E_k$}}
  \subsubsection{The case $M=0$}
  We start with $M=0$. Then, we have
  $$\hat{h}_{1,k}^{(0)} = \hat{h}_{1,(k,0)},~ \hat{h}_{2,k}^{(0)} = \hat{h}_{2,(k,0)}, ~ P_k^{(0 \times 0)}=1, ~k\geq 0, $$
  $$ \frac{d}{dt} \hat{h}_{1,k}^{(0)} = \frac{d}{dt} \hat{h}_{1,(k,0)} = 0, ~ \frac{d}{dt} \hat{h}_{2,k}^{(0)} = \frac{d}{dt} \hat{h}_{2,(k,0)} = 0,   $$
which we get from \eqref{BGKend} with $$L_{1,1}^{(0 \times 0)} = L_{2,1}^{(0 \times 0)} = L_{1,2}^{(0 \times 0)} = L_{2,2}^{(0 \times 0)} = L_{1,3}^{(0 \times 0)} = L_{1,4}^{(0 \times 0)} = L_{2,3}^{(0 \times 0)} = L_{2,4}^{(0 \times 0)}=0. $$

   \subsubsection{The case $M=1$}
   Next, we consider $M=1$. Note that for different M's the time evolution of the common components are different. So for example in the case $M=0$ the time evolution of $\hat{h}_{1,(k,0)}$ is different than the time evolution of $\hat{h}_{1,(k,0)}$ in the case $M=1$. So actually, the two functions are different and also the components should have an index $M$. But in order to not make the notation completely full of indices, we omit this index in the components and use the same notation for all $M$. We will never mix the different cases and then the notation is simpler. In the case $M=1$, we have
  $$\hat{h}_{1,k}^{(1)} = \begin{pmatrix}
  \hat{h}_{1,(k,0)} \\ \hat{h}_{1,(k,1)}
  \end{pmatrix}, ~ \hat{h}_{2,k}^{(1)} = \begin{pmatrix}
  \hat{h}_{2,(k,0)} \\ \hat{h}_{2,(k,1)}
  \end{pmatrix},  ~ k \geq 0, $$
  $$ P_0^{(1 \times 1)} = \textbf{1}_{2 \times 2}, \, P_k^{(1 \times 1)} = \begin{pmatrix}
  1 & - \frac{i \alpha}{k} \\ i \frac{\alpha}{k} & 1
  \end{pmatrix}, ~ k>0. $$
  With equation \eqref{BGKend}, we obtain
 \begin{align*}
 \begin{split}
  &\frac{d}{dt} \hat{h}_{1,k}^{(1)} = \frac{d}{dt}\begin{pmatrix}
  \hat{h}_{1,(k,0)} \\ \hat{h}_{1,(k,1)}
  \end{pmatrix} = - ik \frac{2 \pi}{L} \begin{pmatrix}
  0 & 1 \\ 1 & 0
  \end{pmatrix} \begin{pmatrix}
  \hat{h}_{1,(k,0)} \\ \hat{h}_{1,(k,1)}
  \end{pmatrix} - \nu_{12} n_{\infty, 2} \begin{pmatrix}
  0 & 0 \\ 0 &(1- \delta)
\end{pmatrix}\begin{pmatrix}
  \hat{h}_{1,(k,0)} \\ \hat{h}_{1,(k,1)}
  \end{pmatrix} \\&\hspace{4cm} + \nu_{12} n_{\infty,2} \begin{pmatrix}
  0 & 0 \\ 0 & (1- \delta)\frac{n_{\infty,1}}{n_{\infty,2}} \frac{\sqrt{m_1}}{\sqrt{m_2}} 
\end{pmatrix}\begin{pmatrix}
  \hat{h}_{2,(k,0)} \\ \hat{h}_{2,(k,1)}
  \end{pmatrix},
  \end{split}
  \end{align*}
and 
   \begin{align*}
\begin{split}
 & \frac{d}{dt} \hat{h}_{2,k}^{(1)} = \frac{d}{dt}\begin{pmatrix}
  \hat{h}_{2,(k,0)} \\ \hat{h}_{2,(k,1)}
  \end{pmatrix} = - ik \frac{2 \pi}{L} \begin{pmatrix}
  0 & 1 \\ 1 & 0
  \end{pmatrix} \begin{pmatrix}
  \hat{h}_{2,(k,0)} \\ \hat{h}_{2,(k,1)}
  \end{pmatrix} - \nu_{21} n_{\infty, 1} \begin{pmatrix}
  0 & 0 \\ 0 & \frac{m_1}{m_2} \varepsilon (1- \delta)
\end{pmatrix}\begin{pmatrix}
  \hat{h}_{2,(k,0)} \\ \hat{h}_{2,(k,1)}
  \end{pmatrix} \\&\hspace{4cm}+ \nu_{21} n_{\infty,1} \begin{pmatrix}
  0 & 0 \\ 0 & \varepsilon (1- \delta)\frac{n_{\infty,2}}{n_{\infty,1}} \frac{\sqrt{m_1}}{\sqrt{m_2}} 
\end{pmatrix}\begin{pmatrix}
  \hat{h}_{1,(k,0)} \\ \hat{h}_{1,(k,1)}
  \end{pmatrix}. 
  \end{split}
  \end{align*}
Therefore, we have for $k>0$, 
   \begin{align}
  \begin{split}
 &\textstyle\quad \frac{d}{dt}\hat{h}_{1,k}^{(1)} \cdot P_k^{(1 \times 1)} \cdot \hat{h}_{1,k}^{(1)} +  \hat{h}_{1,k}^{(1)} \cdot P_k^{(1 \times 1)} \cdot \frac{d}{dt}\hat{h}_{1,k}^{(1)} \\& \scriptstyle  \begin{pmatrix}
\scriptstyle    ik \frac{2 \pi}{L} \overline{\hat{h}_{1,(k,1)}} & \scriptstyle i k \frac{2 \pi}{L} \overline{\hat{h}_{1,(k,0)}} \end{pmatrix} \begin{pmatrix}
\scriptstyle 1 & \scriptstyle - \frac{i \tilde{\alpha}}{k} \\ \scriptstyle \frac{i \tilde{\alpha}}{k} & \scriptstyle 1 \end{pmatrix} \begin{pmatrix}
\scriptstyle \hat{h}_{1,(k,0)} \\ \scriptstyle \hat{h}_{1,(k,1)} 
\end{pmatrix} + \begin{pmatrix} \scriptstyle
\overline{\hat{h}_{1,(k,0)}} & \scriptstyle \overline{\hat{h}_{1,(k,1)}} \end{pmatrix} \begin{pmatrix}
\scriptstyle 1 & \scriptstyle - \frac{i \tilde{\alpha}}{k} \\ \scriptstyle \frac{i \tilde{\alpha}}{k} & \scriptstyle 1 \end{pmatrix} \begin{pmatrix}
\scriptstyle - ik \frac{2 \pi}{L} \hat{h}_{1,(k,1)} \\ \scriptstyle -ik \frac{2 \pi}{L} \hat{h}_{1,(k,0)} 
\end{pmatrix} \\ &\scriptstyle + \begin{pmatrix} \scriptstyle 0 & \scriptstyle - \nu_{12} n_{\infty,2} (1- \delta) \overline{\hat{h}_{1,(k,1)}} 
\end{pmatrix} \begin{pmatrix}
\scriptstyle 1 & \scriptstyle - \frac{i \tilde{\alpha}}{k} \\ \scriptstyle \frac{i \tilde{\alpha}}{k} & \scriptstyle 1 \end{pmatrix} \begin{pmatrix}
\scriptstyle \hat{h}_{1,(k,0)} \\ \scriptstyle \hat{h}_{1,(k,1)} 
\end{pmatrix} + \begin{pmatrix}
\scriptstyle     \overline{\hat{h}_{1,(k,0)}} & \scriptstyle  \overline{\hat{h}_{1,(k,1)}} \end{pmatrix} \begin{pmatrix}
\scriptstyle 1 & \scriptstyle - \frac{i \tilde{\alpha}}{k} \\ \scriptstyle \frac{i \tilde{\alpha}}{k} & \scriptstyle 1 \end{pmatrix} \begin{pmatrix} \scriptstyle 0 \\ \scriptstyle - \nu_{12} n_{\infty,2} (1- \delta) \hat{h}_{1,(k,1)} \end{pmatrix} 
\\ &\scriptstyle + \begin{pmatrix} \scriptstyle 0 & \scriptstyle  \nu_{12} n_{\infty,1} (1- \delta)\frac{\sqrt{m_1}}{\sqrt{m_2}} \overline{\hat{h}_{2,(k,1)}} 
\end{pmatrix} \begin{pmatrix}
\scriptstyle 1 & \scriptstyle - \frac{i \tilde{\alpha}}{k} \\ \scriptstyle \frac{i \tilde{\alpha}}{k} & \scriptstyle 1 \end{pmatrix} \begin{pmatrix}
\scriptstyle \hat{h}_{1,(k,0)} \\ \scriptstyle \hat{h}_{1,(k,1)} 
\end{pmatrix} + \begin{pmatrix}
\scriptstyle     \overline{\hat{h}_{1,(k,0)}} & \scriptstyle  \overline{\hat{h}_{1,(k,1)}} \end{pmatrix} \begin{pmatrix}
\scriptstyle 1 & \scriptstyle - \frac{i \tilde{\alpha}}{k} \\ \scriptstyle \frac{i \tilde{\alpha}}{k} & \scriptstyle 1 \end{pmatrix} \begin{pmatrix} \scriptstyle 0 \\ \scriptstyle  \nu_{12} n_{\infty,1} (1- \delta) \frac{\sqrt{m_1}}{\sqrt{m_2}} \hat{h}_{2,(k,1)} \end{pmatrix} 
  \end{split} 
  \label{eq1}
  \end{align}
and {\small\begin{align}
  \begin{split}
 &\textstyle\quad \frac{d}{dt}\hat{h}_{2,k}^{(1)} \cdot P_k^{(1 \times 1)} \cdot \hat{h}_{2,k}^{(1)} +  \hat{h}_{2,k}^{(1)} \cdot P_k^{(1 \times 1)} \cdot \frac{d}{dt}\hat{h}_{2,k}^{(1)} \\& \scriptstyle  \begin{pmatrix}
\scriptstyle    ik \frac{2 \pi}{L} \overline{\hat{h}_{2,(k,1)}} & \scriptstyle i k \frac{2 \pi}{L} \overline{\hat{h}_{2,(k,0)}} \end{pmatrix} \begin{pmatrix}
\scriptstyle 1 & \scriptstyle - \frac{i \tilde{\alpha}}{k} \\ \scriptstyle \frac{i \tilde{\alpha}}{k} & \scriptstyle 1 \end{pmatrix} \begin{pmatrix}
\scriptstyle \hat{h}_{2,(k,0)} \\ \scriptstyle \hat{h}_{2,(k,1)} 
\end{pmatrix} + \begin{pmatrix} \scriptstyle
\overline{\hat{h}_{2,(k,0)}} & \scriptstyle \overline{\hat{h}_{2,(k,1)}} \end{pmatrix} \begin{pmatrix}
\scriptstyle 1 & \scriptstyle - \frac{i \tilde{\alpha}}{k} \\ \scriptstyle \frac{i \tilde{\alpha}}{k} & \scriptstyle 1 \end{pmatrix} \begin{pmatrix}
\scriptstyle - ik \frac{2 \pi}{L} \hat{h}_{2,(k,1)} \\ \scriptstyle -ik \frac{2 \pi}{L} \hat{h}_{2,(k,0)} 
\end{pmatrix} \\ &\scriptstyle + \begin{pmatrix} \scriptstyle 0 & \scriptstyle - \nu_{12} n_{\infty,1} (1- \delta)\frac{m_1}{m_2} \overline{\hat{h}_{2,(k,1)}} 
\end{pmatrix} \begin{pmatrix}
\scriptstyle 1 & \scriptstyle - \frac{i \tilde{\alpha}}{k} \\ \scriptstyle \frac{i \tilde{\alpha}}{k} & \scriptstyle 1 \end{pmatrix} \begin{pmatrix}
\scriptstyle \hat{h}_{2,(k,0)} \\ \scriptstyle \hat{h}_{2,(k,1)} 
\end{pmatrix} + \begin{pmatrix}
\scriptstyle     \overline{\hat{h}_{2,(k,0)}} & \scriptstyle  \overline{\hat{h}_{2,(k,1)}} \end{pmatrix} \begin{pmatrix}
\scriptstyle 1 & \scriptstyle - \frac{i \tilde{\alpha}}{k} \\ \scriptstyle \frac{i \tilde{\alpha}}{k} & \scriptstyle 1 \end{pmatrix} \begin{pmatrix} \scriptstyle 0 \\ \scriptstyle - \nu_{12} n_{\infty,1} (1- \delta) \frac{m_1}{m_2} \hat{h}_{2,(k,1)} \end{pmatrix} 
\\ &\scriptstyle + \begin{pmatrix} \scriptstyle 0 & \scriptstyle  \nu_{12} n_{\infty,2} (1- \delta) \frac{\sqrt{m_1}}{\sqrt{m_2}} \overline{\hat{h}_{1,(k,1)}} 
\end{pmatrix} \begin{pmatrix}
\scriptstyle 1 & \scriptstyle - \frac{i \tilde{\alpha}}{k} \\ \scriptstyle \frac{i \tilde{\alpha}}{k} & \scriptstyle 1 \end{pmatrix} \begin{pmatrix}
\scriptstyle \hat{h}_{2,(k,0)} \\ \scriptstyle \hat{h}_{2,(k,1)} 
\end{pmatrix} + \begin{pmatrix}
\scriptstyle     \overline{\hat{h}_{2,(k,0)}} & \scriptstyle  \overline{\hat{h}_{2,(k,1)}} \end{pmatrix} \begin{pmatrix}
\scriptstyle 1 & \scriptstyle - \frac{i \tilde{\alpha}}{k} \\ \scriptstyle \frac{i \tilde{\alpha}}{k} & \scriptstyle 1 \end{pmatrix} \begin{pmatrix} \scriptstyle 0 \\ \scriptstyle  \nu_{12} n_{\infty,2} (1- \delta) \frac{\sqrt{m_1}}{\sqrt{m_2}} \hat{h}_{1,(k,1)} \end{pmatrix} 
  \end{split} 
  \label{eq2}
  \end{align}}

For the second species, we used \eqref{coll}. We first consider the new terms the new terms compared to the one species case, the third and fourth line of \eqref{eq1} and the third and fourth line of \eqref{eq2}.
  
If we compute the two new lines on the right-hand side of \eqref{eq1}, we get
   \begin{align}
   \begin{split}
& \nu_{12} n_{\infty,2} (1- \delta) \frac{i \tilde{\alpha}}{k} ( \hat{h}_{1,(k,1)} \overline{\hat{h}_{1,(k,0)}} - \hat{h}_{1,(k,0)} \overline{\hat{h}_{1,(k,1)}} ) - 2 \nu_{12} n_{\infty,2} (1- \delta) |\hat{h}_{1,(k,1)}|^2 \\&+ \nu_{12} n_{\infty,1} (1-\delta) \frac{\sqrt{m_1}}{\sqrt{m_2}} i \frac{\tilde{\alpha}}{k} ( \hat{h}_{1,(k,0)} \overline{\hat{h}_{2,(k,1)}} - \hat{h}_{2,(k,1)} \overline{\hat{h}_{1,(k,0)}} ) \\&+ \nu_{12} n_{\infty,1} (1- \delta) \frac{\sqrt{m_1}}{\sqrt{m_2}} ( \hat{h}_{2,(k,1)} \overline{\hat{h}_{1,(k,1)}} + \hat{h}_{1,(k,1)}  \overline{\hat{h}_{2,(k,1)}}) \\ &= - \nu_{12} n_{\infty,2} (1- \delta) 2 \frac{\tilde{\alpha}}{k} \text{Im}( \hat{h}_{1,(k,1)} \overline{\hat{h}_{1,(k,0)}}) - 2 \nu_{12} n_{\infty,2} (1- \delta) |\hat{h}_{1,(k,1)}|^2 \\&- \nu_{12} n_{\infty,1} (1- \delta) \frac{\sqrt{m_1}}{\sqrt{m_2}} 2 \frac{\tilde{\alpha}}{k} \text{Im}(\hat{h}_{1,(k,0)} \overline{\hat{h}_{2,(k,1)}}) + \nu_{12} n_{\infty,1} (1- \delta) \frac{\sqrt{m_1}}{\sqrt{m_2}} 2 \text{Re}(\hat{h}_{2,(k,1)} \overline{\hat{h}_{1,(k,1)}} )
   \end{split}
   \label{eq3}
   \end{align}
  If we compute the two new lines on the right-hand side of \eqref{eq2}, we get 
  \begin{align}
  \begin{split}
- \nu_{12} n_{\infty,1} \frac{m_1}{m_2} (1- \delta) 2 \frac{\tilde{\alpha}}{k} \text{Im}( \overline{\hat{h}_{2,(k,0)}} \hat{h}_{2,(k,1)}) - 2 \nu_{12} n_{\infty,1} \frac{m_1}{m_2} (1- \delta) | \hat{h}_{2,(k,1)}|^2 \\- \nu_{12} n_{\infty,2} \frac{\sqrt{m_1}}{\sqrt{m_2}} \frac{\tilde{\alpha}}{k} (1- \delta) \text{Im} ( \overline{\hat{h}_{1,(k,1)}} \hat{h}_{2,(k,0)} ) + 2 \nu_{12} n_{\infty,2} (1- \delta) \frac{\sqrt{m_1}}{\sqrt{m_2}} \text{Re}( \hat{h}_{1,(k,1)} \overline{\hat{h}_{2,(k,1)}})
  \end{split}
  \label{eq4}
  \end{align}
 Now, we multiply \eqref{eq3} by $\frac{1}{n_{\infty,1}}$ and \eqref{eq4} by $\frac{1}{n_{\infty,2}}$ and add the resulting terms.  Then  we obtain
{\small \begin{align}
 \begin{split}
 - 2 \nu_{12} \frac{n_{\infty,2}}{n_{\infty,1}}  (1- \delta) |\hat{h}_{1,(k,1)}|^2 - 2 \nu_{12} \frac{n_{\infty,1}}{n_{\infty,2}} (1- \delta) \frac{m_1}{m_2} |\hat{h}_{2,(k,1)}|^2 + 4 \nu_{12} (1- \delta) \frac{\sqrt{m_1}}{\sqrt{m_2}} \text{Re}(\hat{h}_{2,(k,1)} \overline{\hat{h}_{1,(k,1)}}) \\ - \nu_{12} (1- \delta) 2 \frac{\tilde{\alpha}}{k} \text{Im}\Big( \frac{n_{\infty,1}}{n_{\infty,2}} \frac{m_1}{m_2} \overline{\hat{h}_{2,(k,0)}} \hat{h}_{2,(k,1)} + \frac{\sqrt{m_1}}{\sqrt{m_2}} \overline{\hat{h}_{1,(k,1)}} \hat{h}_{2,(k,0)} + \frac{n_{\infty,2}}{n_{\infty,1}} \hat{h}_{1,(k,1)} \overline{\hat{h}_{1,(k,0)}} + \frac{\sqrt{m_1}}{\sqrt{m_2}} \hat{h}_{1,(k,0)} \overline{\hat{h}_{2,(k,1)}}\Big)
   \end{split}
   \label{eq5}
\end{align}  }
Now, we use that $\Big|\sqrt{\frac{n_{\infty,2}}{n_{\infty,1}}} \frac{\hat{h}_{1,(k,1)}}{\sqrt{m_1}} -  \sqrt{\frac{n_{\infty,1}}{n_{\infty,2}}}\frac{\hat{h}_{2,(k,1)}}{\sqrt{m_2}}\Big|^2 \geq 0$ and obtain
  $$2   \frac{1}{\sqrt{m_1}\sqrt{m_2}} \text{Re}\Big( \overline{\hat{h}_{1,(k,1)}} \hat{h}_{2,(k,1)} \Big) \leq
  \frac{n_{\infty,2}}{n_{\infty,1}} \frac{|\hat{h}_{1,(k,1)}|^2}{m_1} +  \frac{n_{\infty,1}}{n_{\infty,2}}\frac{|\hat{h}_{2,(k,1)}|^2}{m_2}. $$
With this estimate the first line of the right-hand side of \eqref{eq5} can be bounded above by $0$.


 \textcolor{black}{In conclusion, we get that $\frac{1}{n_{\infty,1}} \eqref{eq1}+\frac{1}{n_{\infty,2}} \eqref{eq2}$ can be bounded from above by the one species terms and the additional imaginary terms coming from  \eqref{eq5}. The expressions \eqref{eq3} and \eqref{eq4} are part of the right-hand side of equation \eqref{eq1} and \eqref{eq2}. Since we multiplied \eqref{eq3} by the weight $\frac{1}{n_{\infty,1}}$ and \eqref{eq4} by $\frac{1}{n_{\infty,2}}$ and summed them up in order to get an estimate from above, we also have to do this for the left-hand side of \eqref{eq1} and \eqref{eq2}, respectively. Therefore, we obtain that}
  $$\frac{1}{n_{\infty,1}} \left(\frac{d}{dt} \hat{h}_{1,k}^{(1)} \cdot P_k^{(1 \times 1)} \cdot \hat{h}_{1,k}^{(1)}+\hat{h}_{1,k}^{(1)} \cdot P_k^{(1 \times 1)} \cdot \frac{d}{dt} \hat{h}_{1,k}^{(1)}\right)+\frac{1}{n_{\infty,2}}\left(\frac{d}{dt} \hat{h}_{2,k}^{(1)} \cdot P_k^{(1 \times 1)} \cdot \hat{h}_{2,k}^{(1)}+\hat{h}_{2,k}^{(1)} \cdot P_k^{(1 \times 1)} \cdot \frac{d}{dt} \hat{h}_{2,k}^{(1)}\right)$$
 can be bounded above by the one species terms and the additional imaginary terms coming from  \eqref{eq5}. \textcolor{black}{Note, that multiplication by $\frac{1}{n_{\infty,1}}$ and $\frac{1}{n_{\infty,2}}$ is not the only possibility to obtain the estimate from above. Another possibility would be to multiply \eqref{eq3} by $\frac{n_{\infty,2}}{n_{\infty,1}+n_{\infty,2}}$ and \eqref{eq4} by $\frac{n_{\infty,1}}{n_{\infty,1}+n_{\infty,2}}$. Then, this would lead to the entropy described in remark \ref{rem}.}
 
  
\subsubsection{The case $M=2$}
  Next, consider $M=2$. Then, we have
  $$\hat{h}_{1,k}^{(2)} = \begin{pmatrix}
  \hat{h}_{1,(k,0)} \\ \hat{h}_{1,(k,1)} \\ \hat{h}_{1,(k,2)}
  \end{pmatrix}, ~ \hat{h}_{2,k}^{(2)} = \begin{pmatrix}
  \hat{h}_{2,(k,0)} \\ \hat{h}_{2,(k,1)} \\ \hat{h}_{2,(k,2)}
  \end{pmatrix}, ~ k \geq 0$$
  $$ P_0^{(2 \times 2)} = \textbf{1}_{3 \times 3},~ P_k^{(2 \times 2)} = \begin{pmatrix}
  1 & - \frac{i \tilde{\alpha}}{k} & 0,  \\ i \frac{\tilde{\alpha}}{k} & 1 & - \frac{i \beta}{k} \\ 0 & \frac{i \beta}{k} & 1
  \end{pmatrix}, k>0. $$
  With equation \eqref{BGKend}, we obtain
   \begin{align*}
   \begin{split}
& \frac{d}{dt} \hat{h}_{1,k}^{(2)} = \frac{d}{dt}\begin{pmatrix}
  \hat{h}_{1,(k,0)} \\ \hat{h}_{1,(k,1)} \\ \hat{h}_{1,(k,2)}
  \end{pmatrix} = - ik \frac{2 \pi}{L} \begin{pmatrix}
  0 & 1 & 0\\ 1 & 0 & \sqrt{2} \\ 0 & \sqrt{2}& 0
  \end{pmatrix} \begin{pmatrix}
  \hat{h}_{1,(k,0)} \\ \hat{h}_{1,(k,1)} \\ \hat{h}_{1,(k,2)}
  \end{pmatrix} \\&\hspace{4cm}- \nu_{12} n_{\infty, 2} \begin{pmatrix}
  0 & 0 & 0 \\ 0 &(1- \delta) & 0 \\ 0 & 0 & (1-\alpha)
\end{pmatrix}\begin{pmatrix}
  \hat{h}_{1,(k,0)} \\ \hat{h}_{1,(k,1)} \\ \hat{h}_{1,(k,2)}
  \end{pmatrix} 
  \\&\hspace{4cm} + \nu_{12} n_{\infty,2} \begin{pmatrix}
  0 & 0 & 0\\ 0 & (1- \delta)\frac{n_{\infty,1}}{n_{\infty,2}} \frac{\sqrt{m_1}}{\sqrt{m_2}} & 0 \\ 0 & 0 & \frac{n_{\infty,1}}{n_{\infty,2}} (1- \alpha)
\end{pmatrix}\begin{pmatrix}
  \hat{h}_{2,(k,0)} \\ \hat{h}_{2,(k,1)} \\ \hat{h}_{2,(k,2)}
  \end{pmatrix},
\end{split}
\end{align*}
and
\begin{align*}
\begin{split}
 & \frac{d}{dt} \hat{h}_{2,k}^{(1)} = \frac{d}{dt}\begin{pmatrix}
  \hat{h}_{2,(k,0)} \\ \hat{h}_{2,(k,1)} \\ \hat{h}_{2,(k,2)}
  \end{pmatrix} = - ik \frac{2 \pi}{L} \begin{pmatrix}
  0 & 1 & 0 \\ 1 & 0 & \sqrt{2} \\ 0 & \sqrt{2} & 0
  \end{pmatrix} \begin{pmatrix}
  \hat{h}_{2,(k,0)} \\ \hat{h}_{2,(k,1)} \\ \hat{h}_{2,(k,2)}
  \end{pmatrix} \\&\hspace{4cm}- \nu_{21} n_{\infty, 1} \begin{pmatrix}
  0 & 0 & 0 \\ 0 & \frac{m_1}{m_2} \varepsilon (1- \delta) & 0 \\ 0 & 0 & \varepsilon (1 - \alpha)
\end{pmatrix}\begin{pmatrix}
  \hat{h}_{2,(k,0)} \\ \hat{h}_{2,(k,1)} \\ \hat{h}_{2,(k,2)}
  \end{pmatrix} \\& \hspace{4cm} + \nu_{21} n_{\infty,1} \begin{pmatrix}
  0 & 0 & 0 \\ 0 & \varepsilon (1- \delta)\frac{n_{\infty,2}}{n_{\infty,1}} \frac{\sqrt{m_1}}{\sqrt{m_2}}  & 0 \\ 0 & 0 & \frac{n_{\infty,2}}{n_{\infty,1}}\varepsilon (1- \alpha)
\end{pmatrix}\begin{pmatrix}
  \hat{h}_{1,(k,0)} \\ \hat{h}_{1,(k,1)} \\ \hat{h}_{1,(k,2)}
  \end{pmatrix}. 
  \end{split}
  \end{align*}
  Therefore, we get for $k>0$
  {\small
    \begin{align}
  \begin{split}
&\quad\frac{d}{dt}\hat{h}_{1,k}^{(2)} \cdot P_k^{(2 \times 2)} \cdot\hat{h}_{1,k}^{(2)} +  \hat{h}_{1,k}^{(2)} \cdot P_k^{(2 \times 2)} \cdot \frac{d}{dt}\hat{h}_{1,k}^{(2)}\\ &= \begin{pmatrix}
i k \frac{2 \pi}{L} \overline{\hat{h}_{1,(k,1)}} & ik \frac{2 \pi}{L} \overline{\hat{h}_{1,(k,0)}} + \sqrt{2 } i k \frac{2 \pi}{L} \overline{\hat{h}_{1,(k,2)}} & i k \frac{2 \pi}{L}\sqrt{2} \overline{\hat{h}_{1,(k,1)}} \end{pmatrix} \begin{pmatrix}
 1 & - \frac{i \tilde{\alpha}}{k} & 0 \\ \frac{i \tilde{\alpha}}{k} & 1 & - \frac{i \beta}{k} \\ 0 & \frac{i \beta}{k} & 1
\end{pmatrix} \begin{pmatrix}
\hat{h}_{1,(k,0)} \\ \hat{h}_{1,(k,1)} \\ \hat{h}_{1,(k,2)}
\end{pmatrix} \\&+ \begin{pmatrix}
\overline{\hat{h}_{1,(k,0)}} & \overline{\hat{h}_{1,(k,1)}} & \overline{\hat{h}_{1,(k,2)}} 
\end{pmatrix} \begin{pmatrix}
 1 & - \frac{i \tilde{\alpha}}{k} & 0 \\ \frac{i \tilde{\alpha}}{k} & 1 & - \frac{i \beta}{k} \\ 0 & \frac{i \beta}{k} & 1
\end{pmatrix}  \begin{pmatrix}
-i k \frac{2 \pi}{L} \hat{h}_{1,(k,1)} \\- ik \frac{2 \pi}{L} \hat{h}_{1,(k,0)} - \sqrt{2 } i k \frac{2 \pi}{L} \hat{h}_{1,(k,2)} \\ -i k \frac{2 \pi}{L}\sqrt{2} \hat{h}_{1,(k,1)} \end{pmatrix} \\&+ \begin{pmatrix}
0 & - \nu_{12} n_{\infty,2} (1- \delta) \overline{\hat{h}_{1,(k,1)}} & - \nu_{12} n_{\infty,2} (1- \alpha) \overline{\hat{h}_{1,(k,2)}} \end{pmatrix} \begin{pmatrix}
 1 & - \frac{i \tilde{\alpha}}{k} & 0 \\ \frac{i \tilde{\alpha}}{k} & 1 & - \frac{i \beta}{k} \\ 0 & \frac{i \beta}{k} & 1
\end{pmatrix} \begin{pmatrix}
\hat{h}_{1,(k,0)} \\ \hat{h}_{1,(k,1)} \\ \hat{h}_{1,(k,2)}
\end{pmatrix} \\&+ \begin{pmatrix}
\overline{\hat{h}_{1,(k,0)}} & \overline{\hat{h}_{1,(k,1)}} & \overline{\hat{h}_{1,(k,2)}} 
\end{pmatrix} \begin{pmatrix}
 1 & - \frac{i \tilde{\alpha}}{k} & 0 \\ \frac{i \tilde{\alpha}}{k} & 1 & - \frac{i \beta}{k} \\ 0 & \frac{i \beta}{k} & 1
\end{pmatrix}  \begin{pmatrix}
0 \\- \nu_{12} n_{\infty,2} (1- \delta) \hat{h}_{1,(k,1)}\\ - \nu_{12} n_{\infty,2} (1- \alpha) \hat{h}_{1,(k,2)} \end{pmatrix} \\ &+
\begin{pmatrix}
0 & \nu_{12} n_{\infty,1} (1- \delta) \frac{\sqrt{m_1}}{\sqrt{m_2}} \overline{\hat{h}_{2,(k,1)}} & \nu_{12} n_{\infty,1} (1- \alpha) \overline{\hat{h}_{2,(k,2)}} \end{pmatrix} \begin{pmatrix}
 1 & - \frac{i \tilde{\alpha}}{k} & 0 \\ \frac{i \tilde{\alpha}}{k} & 1 & - \frac{i \beta}{k} \\ 0 & \frac{i \beta}{k} & 1
\end{pmatrix} \begin{pmatrix}
\hat{h}_{1,(k,0)} \\ \hat{h}_{1,(k,1)} \\ \hat{h}_{1,(k,2)}
\end{pmatrix} \\&+ \begin{pmatrix}
\overline{\hat{h}_{1,(k,0)}} & \overline{\hat{h}_{1,(k,1)}} & \overline{\hat{h}_{1,(k,2)}} 
\end{pmatrix} \begin{pmatrix}
 1 & - \frac{i \tilde{\alpha}}{k} & 0 \\ \frac{i \tilde{\alpha}}{k} & 1 & - \frac{i \beta}{k} \\ 0 & \frac{i \beta}{k} & 1
\end{pmatrix}  \begin{pmatrix}
0 \\\nu_{12} n_{\infty,1} (1- \delta) \frac{\sqrt{m_1}}{\sqrt{m_2}} \hat{h}_{2,(k,1)} \\ \nu_{12} n_{\infty,1} (1- \alpha) \hat{h}_{2,(k,2)} \end{pmatrix}
  \end{split} 
  \label{eq6}
  \end{align}}
and
  {\small
     \begin{align}
  \begin{split}
&\quad\frac{d}{dt}\hat{h}_{2,k}^{(2)} \cdot P_k^{(2 \times 2)} \cdot \hat{h}_{2,k}^{(2)} +  \hat{h}_{2,k}^{(2)} \cdot P_k^{(2 \times 2)} \cdot \frac{d}{dt}\hat{h}_{2,k}^{(2)} \\&= \begin{pmatrix}
i k \frac{2 \pi}{L} \overline{\hat{h}_{2,(k,1)}} & ik \frac{2 \pi}{L} \overline{\hat{h}_{2,(k,0)}} + \sqrt{2 } i k \frac{2 \pi}{L} \overline{\hat{h}_{2,(k,2)}} & i k \frac{2 \pi}{L}\sqrt{2} \overline{\hat{h}_{2,(k,1)}} \end{pmatrix} \begin{pmatrix}
 1 & - \frac{i \tilde{\alpha}}{k} & 0 \\ \frac{i \tilde{\alpha}}{k} & 1 & - \frac{i \beta}{k} \\ 0 & \frac{i \beta}{k} & 1
\end{pmatrix} \begin{pmatrix}
\hat{h}_{2,(k,0)} \\ \hat{h}_{2,(k,1)} \\ \hat{h}_{2,(k,2)}
\end{pmatrix} \\&+ \begin{pmatrix}
\overline{\hat{h}_{2,(k,0)}} & \overline{\hat{h}_{2,(k,1)}} & \overline{\hat{h}_{2,(k,2)}} 
\end{pmatrix} \begin{pmatrix}
 1 & - \frac{i \tilde{\alpha}}{k} & 0 \\ \frac{i \tilde{\alpha}}{k} & 1 & - \frac{i \beta}{k} \\ 0 & \frac{i \beta}{k} & 1
\end{pmatrix}  \begin{pmatrix}
-i k \frac{2 \pi}{L} \hat{h}_{2,(k,1)} \\- ik \frac{2 \pi}{L} \hat{h}_{2,(k,0)} - \sqrt{2 } i k \frac{2 \pi}{L} \hat{h}_{2,(k,2)} \\ -i k \frac{2 \pi}{L}\sqrt{2} \hat{h}_{2,(k,1)} \end{pmatrix} \\&+ \begin{pmatrix}
0 & - \nu_{12} n_{\infty,1} \frac{m_1}{m_2}(1- \delta) \overline{\hat{h}_{2,(k,1)}} & - \nu_{12} n_{\infty,1} (1- \alpha) \overline{\hat{h}_{2,(k,2)}} \end{pmatrix} \begin{pmatrix}
 1 & - \frac{i \tilde{\alpha}}{k} & 0 \\ \frac{i \tilde{\alpha}}{k} & 1 & - \frac{i \beta}{k} \\ 0 & \frac{i \beta}{k} & 1
\end{pmatrix} \begin{pmatrix}
\hat{h}_{2,(k,0)} \\ \hat{h}_{2,(k,1)} \\ \hat{h}_{2,(k,2)}
\end{pmatrix} \\&+ \begin{pmatrix}
\overline{\hat{h}_{2,(k,0)}} & \overline{\hat{h}_{2,(k,1)}} & \overline{\hat{h}_{2,(k,2)}} 
\end{pmatrix} \begin{pmatrix}
 1 & - \frac{i \tilde{\alpha}}{k} & 0 \\ \frac{i \tilde{\alpha}}{k} & 1 & - \frac{i \beta}{k} \\ 0 & \frac{i \beta}{k} & 1
\end{pmatrix}  \begin{pmatrix}
0 \\- \nu_{12} n_{\infty,1} \frac{m_1}{m_2}(1- \delta) \hat{h}_{2,(k,1)}\\ - \nu_{12} n_{\infty,1} (1- \alpha) \hat{h}_{2,(k,2)} \end{pmatrix} \\ &+
\begin{pmatrix}
0 & \nu_{12} n_{\infty,2} (1- \delta) \frac{\sqrt{m_1}}{\sqrt{m_2}} \overline{\hat{h}_{1,(k,1)}} & \nu_{12} n_{\infty,2} (1- \alpha) \overline{\hat{h}_{1,(k,2)}} \end{pmatrix} \begin{pmatrix}
 1 & - \frac{i \tilde{\alpha}}{k} & 0 \\ \frac{i \tilde{\alpha}}{k} & 1 & - \frac{i \beta}{k} \\ 0 & \frac{i \beta}{k} & 1
\end{pmatrix} \begin{pmatrix}
\hat{h}_{2,(k,0)} \\ \hat{h}_{2,(k,1)} \\ \hat{h}_{2,(k,2)}
\end{pmatrix} \\&+ \begin{pmatrix}
\overline{\hat{h}_{2,(k,0)}} & \overline{\hat{h}_{2,(k,1)}} & \overline{\hat{h}_{2,(k,2)}} 
\end{pmatrix} \begin{pmatrix}
 1 & - \frac{i \tilde{\alpha}}{k} & 0 \\ \frac{i \tilde{\alpha}}{k} & 1 & - \frac{i \beta}{k} \\ 0 & \frac{i \beta}{k} & 1
\end{pmatrix}  \begin{pmatrix}
0 \\\nu_{12} n_{\infty,2} (1- \delta) \frac{\sqrt{m_1}}{\sqrt{m_2}} \hat{h}_{1,(k,1)} \\ \nu_{12} n_{\infty,2} (1- \alpha) \hat{h}_{1,(k,2)} \end{pmatrix}
  \end{split} 
  \label{eq7}
  \end{align}}

Again, the first two lines on the right hand side of \eqref{eq6} is treated in the same way as \cite{achleitner2017multi} since this term is the same as in the one species case. 
\textcolor{black}{The same is true} for the first two lines on the right-hand side of \eqref{eq7}. It remains to estimate the remaining terms of \eqref{eq6} and \eqref{eq7}.  
Compute the two remaining terms on the right-hand side of \eqref{eq6}, we get 
\begin{align}
\begin{split}
- 2 \nu_{12} n_{\infty,2} (1- \delta) |\hat{h}_{1,(k,1)}|^2 + \nu_{12} n_{\infty,1} (1- \delta) \frac{\sqrt{m_1}}{\sqrt{m_2}} 2 \text{Re}( \hat{h}_{2,(k,1)} \overline{\hat{h}_{1,(k,1)}} )\\ - 2 \nu_{12} n_{\infty,2} (1- \alpha) |\hat{h}_{1,(k,1)}|^2 + 2 \nu_{12} n_{\infty,1} (1- \alpha) \text{Re}(\overline{\hat{h}_{1,(k,2)}} \hat{h}_{2,(k,2)})\\ - \nu_{12} n_{\infty,2} (1- \delta) 2 \frac{\tilde{\alpha}}{k} \text{Im}(\hat{h}_{1,(k,1)}
 \overline{ \hat{h}_{1,(k,0)}}) - \nu_{12} n_{\infty,1}(1-\delta) \frac{\sqrt{m_1}}{\sqrt{m_2}} 2 \frac{\tilde{\alpha}}{k} \text{Im}(\hat{h}_{1,(k,0)} \overline{\hat{h}_{2,(k,1)}}) \\- 2 \frac{\beta}{k} \nu_{12} n_{\infty,2}  (1- \delta) \text{Im}(\overline{\hat{h}_{1,(k,1)}} \hat{h}_{1,(k,2)}) - 2 \frac{\beta}{k} \nu_{12} n_{\infty,1} (1- \delta) \frac{\sqrt{m_1}}{\sqrt{m_2}} \text{Im}(\overline{\hat{h}_{1,(,k,2)}} \hat{h}_{2,(k,1)})\\ - \frac{\beta}{k} \nu_{12} n_{\infty,2} (1- \alpha) 2 \text{Im}(\overline{\hat{h}_{1,(k,1)}} \hat{h}_{1,(k,2)}) - 2 \frac{\beta}{k} \nu_{12} n_{\infty,1} (1- \alpha) \text{Im}(\overline{\hat{h}_{2,(k,2)}} \hat{h}_{1,(k,1)})\end{split}
\label{eq8}
\end{align}
Compute the two remaining terms on the right-hand side of \eqref{eq7}, we get
\begin{align}
\begin{split}
- 2 \nu_{12} n_{\infty,1} \frac{m_1}{m_2} (1- \delta) |\hat{h}_{2,(k,1)}|^2 + \nu_{12} n_{\infty,2} (1- \delta) \frac{\sqrt{m_1}}{\sqrt{m_2}} 2 \text{Re}( \hat{h}_{1,(k,1)} \overline{\hat{h}_{2,(k,1)}} )\\ - 2 \nu_{12} n_{\infty,1} (1- \alpha) |\hat{h}_{2,(k,2)}|^2 + 2 \nu_{12} n_{\infty,2} (1- \alpha) \text{Re}(\overline{\hat{h}_{1,(k,2)}} \hat{h}_{2,(k,1)})\\ - \nu_{12} n_{\infty,1} (1- \delta) 2 \frac{\tilde{\alpha}}{k} \frac{m_1}{m_2} \text{Im}(\hat{h}_{2,(k,1)}
 \overline{ \hat{h}_{2,(k,0)}}) - \nu_{12} n_{\infty,2}(1-\delta) \frac{\sqrt{m_1}}{\sqrt{m_2}} 2 \frac{\tilde{\alpha}}{k} \text{Im}(\hat{h}_{2,(k,0)} \overline{\hat{h}_{1,(k,1)}}) \\- 2 \frac{\beta}{k} \nu_{12} n_{\infty,1} \frac{m_1}{m_2} (1- \delta) \text{Im}(\overline{\hat{h}_{2,(k,1)}} \hat{h}_{2,(k,2)}) - 2 \frac{\beta}{k} \nu_{12} n_{\infty,2} (1- \delta) \frac{\sqrt{m_1}}{\sqrt{m_2}} \text{Im}(\hat{h}_{1,(k,1)} \overline{\hat{h}_{2,(k,2)}})\\ - \frac{\beta}{k} \nu_{12} n_{\infty,1} (1- \alpha) 2 \text{Im}(\overline{\hat{h}_{2,(k,1)}} \hat{h}_{2,(k,2)}) - 2 \frac{\beta}{k} \nu_{12} n_{\infty,2} (1- \alpha) \text{Im}(\overline{\hat{h}_{2,(k,1)}} \hat{h}_{1,(k,2)})
\end{split}
\label{eq9}
\end{align}
Now, we multiply \eqref{eq8} by $\frac{1}{n_{\infty,1}}$ and \eqref{eq9} by $\frac{1}{n_{\infty,2}}$ and add the resulting terms. In the resulting term, we already estimated the terms with $(1- \delta)$  from above in the case $M=1$. The terms with $(1- \alpha)$ can be estimated by zero from above in the same way as the terms with $(1-\delta)$.

In conclusion, we get that 
$$\frac{1}{n_{\infty,1}}\left(\frac{d}{dt}\hat{h}_{1,k}^{(2)} \cdot P_k^{(2 \times 2)} \cdot \hat{h}_{1,k}^{(2)}+\hat{h}_{1,k}^{(2)} \cdot P_k^{(2 \times 2)} \cdot \frac{d}{dt}\hat{h}_{1,k}^{(2)}\right) + \frac{1}{n_{\infty,2}}\left(\frac{d}{dt} \hat{h}_{2,k}^{(2)} \cdot P_k^{(2 \times 2)} \cdot  \hat{h}_{2,k}^{(2)}+ \hat{h}_{2,k}^{(2)} \cdot P_k^{(2 \times 2)} \cdot \frac{d}{dt} \hat{h}_{2,k}^{(2)}\right)$$
can be bounded from above by the one species terms and by the additional imaginary terms from \eqref{eq8} and \eqref{eq9}.
\subsubsection{The cases $M=3$ and $M>3$}
\label{sec4.1.4}
Next, we consider $M=3$. Then, we have
 $$\hat{h}_{1,k}^{(3)} = \begin{pmatrix}
  \hat{h}_{1,(k,0)} \\ \hat{h}_{1,(k,1)} \\ \hat{h}_{1,(k,2)} \\ \hat{h}_{1,(k,3)}
  \end{pmatrix}, ~ \hat{h}_{2,k}^{(3)} = \begin{pmatrix}
  \hat{h}_{2,(k,0)} \\ \hat{h}_{2,(k,1)} \\ \hat{h}_{2,(k,2)} \\ \hat{h}_{2,(k,3)}
  \end{pmatrix}, ~ k \geq 0, $$
  $$ P_0^{(3 \times 3)} = \textbf{1}_{4 \times 4},~ P_k^{(3 \times 3)} = \begin{pmatrix}
  1 & - \frac{i \tilde{\alpha}}{k} & 0 & 0 \\ i \frac{\tilde{\alpha}}{k} & 1 & - \frac{i \beta}{k} & 0 \\ 0 & \frac{i \beta}{k} & 1 & - \frac{i \tilde{\gamma}}{k} \\ 0 & 0& \frac{i \tilde{\gamma}}{k} & 1
  \end{pmatrix}, k>0. $$
With equation \eqref{BGKend}, we obtain 
\begin{align*}
 &\frac{d}{dt} \hat{h}_{1,k}^{(3)} = \frac{d}{dt}\begin{pmatrix}
  \hat{h}_{1,(k,0)} \\ \hat{h}_{1,(k,1)} \\ \hat{h}_{1,(k,2)} \\ \hat{h}_{1,(k,3)}
  \end{pmatrix} = - ik \frac{2 \pi}{L} \begin{pmatrix}
  0 & 1 & 0 & 0\\ 1 & 0 & \sqrt{2} & 0\\ 0 & \sqrt{2} & 0 & \sqrt{3} \\ 0 & 0& \sqrt{3} & 0
  \end{pmatrix} \begin{pmatrix}
  \hat{h}_{1,(k,0)} \\ \hat{h}_{1,(k,1)} \\ \hat{h}_{1,(k,2)} \\ \hat{h}_{1,(k,3)}
  \end{pmatrix}\\&- \nu_{11} n_{\infty,1} \begin{pmatrix}
  0&0&0&0 \\ 0&0&0&0 \\ 0&0&0&0 \\ 0&0&0&1
\end{pmatrix} \begin{pmatrix}
  \hat{h}_{1,(k,0)} \\ \hat{h}_{1,(k,1)} \\ \hat{h}_{1,(k,2)} \\ \hat{h}_{1,(k,3)}
  \end{pmatrix}
  - \nu_{12} n_{\infty, 2} \begin{pmatrix}
  0 & 0 & 0 & 0\\ 0 &(1- \delta) & 0 & 0\\ 0 & 0 & (1-\alpha) & 0 \\ 0 & 0 & 0 & 1
\end{pmatrix}\begin{pmatrix}
  \hat{h}_{1,(k,0)} \\ \hat{h}_{1,(k,1)} \\ \hat{h}_{1,(k,2)} \\ \hat{h}_{1,(k,3)}
  \end{pmatrix} 
  \\&\hspace{4cm}+ \nu_{12} n_{\infty,2} \begin{pmatrix}
  0 & 0 & 0 & 0\\ 0 & (1- \delta)\frac{n_{\infty,1}}{n_{\infty,2}} \frac{\sqrt{m_1}}{\sqrt{m_2}} & 0 & 0\\ 0 & 0 & \frac{n_{\infty,1}}{n_{\infty,2}} (1- \alpha) & 0 \\ 0 & 0&0&0
\end{pmatrix}\begin{pmatrix}
  \hat{h}_{2,(k,0)} \\ \hat{h}_{2,(k,1)} \\ \hat{h}_{2,(k,2)} \\ \hat{h}_{2,(k,3)}
  \end{pmatrix}, 
  \end{align*}
and
     \begin{align*}
  &\frac{d}{dt} \hat{h}_{2,k}^{(3)} = \frac{d}{dt}\begin{pmatrix}
  \hat{h}_{2,(k,0)} \\ \hat{h}_{2,(k,1)} \\ \hat{h}_{2,(k,2)} \\ \hat{h}_{2,(k,3)}
  \end{pmatrix} = - ik \frac{2 \pi}{L} \begin{pmatrix}
  0 & 1 & 0 & 0 \\ 1 & 0 & \sqrt{2} & 0\\ 0 & \sqrt{2} & 0 & \sqrt{3} \\ 0 & 0 & \sqrt{3} & 0
  \end{pmatrix} \begin{pmatrix}
  \hat{h}_{2,(k,0)} \\ \hat{h}_{2,(k,1)} \\ \hat{h}_{2,(k,2)} \\ \hat{h}_{2,(k,3)}
  \end{pmatrix}\\&- \nu_{22} n_{\infty,2} \begin{pmatrix}
  0&0&0&0 \\ 0&0&0&0 \\ 0&0&0&0 \\ 0&0&0&1
\end{pmatrix} \begin{pmatrix}
  \hat{h}_{2,(k,0)} \\ \hat{h}_{2,(k,1)} \\ \hat{h}_{2,(k,2)} \\ \hat{h}_{2,(k,3)}
  \end{pmatrix}
 - \nu_{21} n_{\infty, 1} \begin{pmatrix}
  0 & 0 & 0 &0\\ 0 & \frac{m_1}{m_2} \varepsilon (1- \delta) & 0 &0\\ 0 & 0 & \varepsilon (1 - \alpha) &0 \\ 0&0&0&1
\end{pmatrix}\begin{pmatrix}
  \hat{h}_{2,(k,0)} \\ \hat{h}_{2,(k,1)} \\ \hat{h}_{2,(k,2)} \\ \hat{h}_{2,(k,3)}
  \end{pmatrix}
  \\&\hspace{4cm} + \nu_{21} n_{\infty,1} \begin{pmatrix}
  0 & 0 & 0 & 0\\ 0 & \varepsilon (1- \delta)\frac{n_{\infty,1}}{n_{\infty,2}} \frac{\sqrt{m_1}}{\sqrt{m_2}}  & 0 & 0\\ 0 & 0 & \varepsilon (1- \alpha) & 0 \\ 0&0&0&0
\end{pmatrix}\begin{pmatrix}
  \hat{h}_{1,(k,0)} \\ \hat{h}_{1,(k,1)} \\ \hat{h}_{1,(k,2)} \\ \hat{h}_{1,(k,3)}
  \end{pmatrix}. 
  \end{align*}
Therefore, we get for species $1$ for $k>0$
{\scriptsize
    \begin{align*}
  \begin{split}
  &\frac{d}{dt}\hat{h}_{1,k}^{(3)} \cdot P_k^{(3 \times 3)} \cdot \hat{h}_{1,k}^{(3)}+\hat{h}_{1,k}^{(3)} \cdot P_k^{(3 \times 3)} \cdot \frac{d}{dt}\hat{h}_{1,k}^{(3)} \\&=  ik \frac{2 \pi}{L}  \begin{pmatrix}
 \overline{ \hat{h}_{1,(k,0)}} & \overline{\hat{h}_{1,(k,1)}} & \overline{\hat{h}_{1,(k,2)}} & \overline{\hat{h}_{1,(k,3)}}
  \end{pmatrix}\begin{pmatrix}
  0 & 1 & 0 & 0\\ 1 & 0 & \sqrt{2} & 0 \\ 0 & \sqrt{2} & 0 & \sqrt{3} \\ 0&0&\sqrt{3}&0
  \end{pmatrix}\begin{pmatrix}
  1 & - \frac{i \tilde{\alpha}}{k} & 0 & 0\\ i \frac{\tilde{\alpha}}{k} & 1 & - \frac{i \beta}{k} & 0\\ 0 & \frac{i \beta}{k} & 1 & - \frac{i \tilde{\gamma}}{k} \\ 0 & 0& \frac{i \tilde{\gamma}}{k} & 1
  \end{pmatrix}\begin{pmatrix} \hat{h}_{1,(k,0)} \\ \hat{h}_{1,(k,1)} \\ \hat{h}_{1,(k,2)} \\ \hat{h}_{1,(k,3)} \end{pmatrix} \\&\quad  \begin{pmatrix} \overline{\hat{h}_{1,(k,0)}} & \overline{\hat{h}_{1,(k,1)}} & \overline{\hat{h}_{1,(k,2)}} & \overline{\hat{h}_{1,(k,3)}} \end{pmatrix} \begin{pmatrix}
  1 & - \frac{i \tilde{\alpha}}{k} & 0 & 0\\ i \frac{\tilde{\alpha}}{k} & 1 & - \frac{i \beta}{k} & 0\\ 0 & \frac{i \beta}{k} & 1 & - \frac{i \tilde{\gamma}}{k} \\ 0 & 0& \frac{i \tilde{\gamma}}{k} & 1
  \end{pmatrix} \begin{pmatrix}
  0 & 1 & 0 & 0\\ 1 & 0 & \sqrt{2} & 0 \\ 0 & \sqrt{2} & 0 & \sqrt{3} \\ 0&0&\sqrt{3}&0
  \end{pmatrix}  \begin{pmatrix}
 - i k \frac{2 \pi}{L}	\hat{h}_{1,(k,0)} \\ - i k \frac{2 \pi}{L} \hat{h}_{1,(k,1)} \\ - i k \frac{2 \pi}{L} \hat{h}_{1,(k,2)} \\ - i k \frac{2 \pi}{L} \hat{h}_{1,(k,3)}
  \end{pmatrix}\\&\quad - \nu_{11} n_{\infty,1} \begin{pmatrix}
 \overline{ \hat{h}_{1,(k,0)}} & \overline{\hat{h}_{1,(k,1)}} & \overline{\hat{h}_{1,(k,2)}} & \overline{\hat{h}_{1,(k,3)}}
  \end{pmatrix} \begin{pmatrix}
0&0&0&0\\0&0&0&0\\0&0&0&0\\0&0&0&1
\end{pmatrix} \begin{pmatrix}
  1 & - \frac{i \tilde{\alpha}}{k} & 0 & 0\\ i \frac{\tilde{\alpha}}{k} & 1 & - \frac{i \beta}{k} & 0\\ 0 & \frac{i \beta}{k} & 1 & - \frac{i \tilde{\gamma}}{k} \\ 0 & 0& \frac{i \tilde{\gamma}}{k} & 1
\end{pmatrix} \begin{pmatrix}
\hat{h}_{1,(k,0)} \\ \hat{h}_{1,(k,1)} \\ \hat{h}_{1,(k,2)} \\ \hat{h}_{1,(k,3)}
\end{pmatrix} \\&\quad - \nu_{11} n_{\infty,1} \begin{pmatrix} \overline{\hat{h}_{1,(k,0)}} & \overline{\hat{h}_{1,(k,1)}} & \overline{\hat{h}_{1,(k,2)}} & \overline{\hat{h}_{1,(k,3)}} \end{pmatrix} \begin{pmatrix}
  1 & - \frac{i \tilde{\alpha}}{k} & 0 & 0\\ i \frac{\tilde{\alpha}}{k} & 1 & - \frac{i \beta}{k} & 0\\ 0 & \frac{i \beta}{k} & 1 & - \frac{i \tilde{\gamma}}{k} \\ 0 & 0& \frac{i \tilde{\gamma}}{k} & 1
\end{pmatrix} \begin{pmatrix}
0&0&0&0\\0&0&0&0\\0&0&0&0\\0&0&0&1
\end{pmatrix}  \begin{pmatrix}
  \hat{h}_{1,(k,0)} \\ \hat{h}_{1,(k,1)} \\ \hat{h}_{1,(k,2)} \\ \hat{h}_{1,(k,3)}
  \end{pmatrix} \\&\quad + \begin{pmatrix} 0 & - \nu_{12} n_{\infty,2} (1- \delta) \overline{\hat{h}_{1,(k,1)}} & - \nu_{12} n_{\infty,2} (1- \alpha) \overline{\hat{h}_{1,(k,2)}} & - \nu_{12} n_{\infty,2} \overline{\hat{h}_{1,(k,3)}} \end{pmatrix} \begin{pmatrix}
  1 & - \frac{i \tilde{\alpha}}{k} & 0 & 0\\ i \frac{\tilde{\alpha}}{k} & 1 & - \frac{i \beta}{k} & 0\\ 0 & \frac{i \beta}{k} & 1 & - \frac{i \tilde{\gamma}}{k} \\ 0 & 0& \frac{i \tilde{\gamma}}{k} & 1
  \end{pmatrix}\begin{pmatrix} \hat{h}_{1,(k,0)} \\ \hat{h}_{1,(k,1)} \\ \hat{h}_{1,(k,2)} \\ \hat{h}_{1,(k,3)} \end{pmatrix} \\ &+ \begin{pmatrix} \overline{\hat{h}_{1,(k,0)}} & \overline{\hat{h}_{1,(k,1)}} & \overline{\hat{h}_{1,(k,2)}} & \overline{\hat{h}_{1,(k,3)}} \end{pmatrix} \begin{pmatrix}
  1 & - \frac{i \tilde{\alpha}}{k} & 0 & 0\\ i \frac{\tilde{\alpha}}{k} & 1 & - \frac{i \beta}{k} & 0\\ 0 & \frac{i \beta}{k} & 1 & - \frac{i \tilde{\gamma}}{k} \\ 0 & 0& \frac{i \tilde{\gamma}}{k} & 1
  \end{pmatrix} \begin{pmatrix}
  0 \\ - \nu_{12} n_{\infty,2} (1- \delta) \hat{h}_{1,(k,1)} \\ - \nu_{12} n_{\infty,2} (1- \alpha) \hat{h}_{1,(k,2)} \\ - \nu_{12} n_{\infty,2} \hat{h}_{1,(k,3)} \end{pmatrix} \\ &+ \begin{pmatrix}
  0 & \nu_{12} n_{\infty,1} (1- \delta) \frac{\sqrt{m_1}}{\sqrt{m_2}} \overline{\hat{h}_{2,(k,1)}} & \nu_{12} n_{\infty,1} (1- \alpha) \overline{\hat{h}_{2,(k,2)}} & 0 
  \end{pmatrix} \begin{pmatrix}
  1 & - \frac{i \tilde{\alpha}}{k} & 0 & 0\\ i \frac{\tilde{\alpha}}{k} & 1 & - \frac{i \beta}{k} & 0\\ 0 & \frac{i \beta}{k} & 1 & - \frac{i \tilde{\gamma}}{k} \\ 0 & 0& \frac{i \tilde{\gamma}}{k} & 1
  \end{pmatrix} \begin{pmatrix}
  \hat{h}_{1,(k,0)} \\ \hat{h}_{1,(k,1)} \\ \hat{h}_{1,(k,2)} \\ \hat{h}_{1,(k,3)}
  \end{pmatrix} \\ &+ \begin{pmatrix}
  \overline{\hat{h}_{1,(k,0)}} & \overline{\hat{h}_{1,(k,1)}} & \overline{\hat{h}_{1,(k,2)}} & \overline{\hat{h}_{1,(k,3)}} 
  \end{pmatrix} \begin{pmatrix}
  1 & - \frac{i \tilde{\alpha}}{k} & 0 & 0\\ i \frac{\tilde{\alpha}}{k} & 1 & - \frac{i \beta}{k} & 0\\ 0 & \frac{i \beta}{k} & 1 & - \frac{i \tilde{\gamma}}{k} \\ 0 & 0& \frac{i \tilde{\gamma}}{k} & 1
  \end{pmatrix} \begin{pmatrix}
  0 \\ \nu_{12} n_{\infty,1} (1- \delta) \frac{\sqrt{m_1}}{\sqrt{m_2}} \hat{h}_{2,(k,1)} \\ \nu_{12} n_{\infty,1} (1- \alpha) \hat{h}_{2,(k,2)} \\ 0
  \end{pmatrix}
  \end{split} 
  \end{align*}}
  The only new terms compared to the one species terms and the case $M=2$ can be computed as
  \begin{align}
  \begin{split}
- \nu_{12} n_{\infty,2} (1- \alpha) \frac{\tilde{\gamma}}{k} 2 \text{Im}(\overline{\hat{h}_{1,(k,2)}} \hat{h}_{1,(k,3)}) - 2 \nu_{12} n_{\infty,2} |\hat{h}_{1,(k,3)}|^2 - \frac{\tilde{\gamma}}{k} \nu_{12} n_{\infty,2} 2 \text{Im}(\overline{\hat{h}_{1,(k,2)}} \hat{h}_{1,(k,3)})\\ -2 \nu_{12} n_{\infty,2} |\hat{h}_{1,(k,3)}|^2 - \frac{\tilde{\gamma}}{k} \nu_{12} (1- \alpha) n_{\infty,1} \text{Im}(\hat{h}_{2,(k,2)} \overline{\hat{h}_{1,(k,3)}})
  \end{split}
  \label{eq10}
  \end{align}
Similar, for the second species. To summarize, we still have to estimate all one species terms and all new imaginary terms coming from \eqref{eq5}, $\frac{1}{n_{\infty,1}} \eqref{eq8} + \frac{1}{n_{\infty,2}} \eqref{eq9}$ and $\frac{1}{n_{\infty,1}}$ times \eqref{eq10} plus $\frac{1}{n_{\infty,2}}$ times the corresponding terms to \eqref{eq10} of species 2. The imaginary terms can be estimated by squared terms by using $- 2 \text{Im}(\overline{z_1} z_2) \leq \frac{1}{c}|z_1|^2 + c |z_2|^2$ for complex numbers $z_1, z_2$ and a real positive constant $c$. This inequality comes from $|\frac{1}{\sqrt{c}} z_1- \sqrt{c} iz_2|^2 \geq 0$. Then, by writing down all one species terms plus the additional terms, one can see that the additional terms can be absorbed into the one species terms  \cite{achleitner2017multi} by a certain choice of the parameters $\tilde{\alpha}, \beta$ and $\tilde{\gamma}$ and assuming that $\nu_{11} n_{\infty,1} + \nu_{12} n_{\infty,2} =1$ and $\nu_{22} n_{\infty,2} + \nu_{21} n_{\infty,1}= 1$.

In conclusion, we get that 
$$\frac{1}{n_{\infty,1}}\left(\frac{d}{dt}\hat{h}_{1,k}^{(3)} \cdot P_k^{(2 \times 2)} \cdot \hat{h}_{1,k}^{(3)}+\hat{h}_{1,k}^{(3)} \cdot P_k^{(3 \times 3)} \cdot \frac{d}{dt}\hat{h}_{1,k}^{(3)}\right) + \frac{1}{n_{\infty,2}}\left(\frac{d}{dt} \hat{h}_{2,k}^{(3)} \cdot P_k^{(3 \times 3)} \cdot  \hat{h}_{2,k}^{(3)}+ \hat{h}_{2,k}^{(3)} \cdot P_k^{(3 \times 3)} \cdot \frac{d}{dt} \hat{h}_{2,k}^{3)}\right)$$
can be bounded from above by 
$$- 2 \mu \left(\frac{1}{n_{\infty,1}} \hat{h}_{1,k}^{(3)} \cdot P_k^{(3 \times 3)} \cdot \hat{h}_{1,k}^{(3)} +\frac{1}{n_{\infty,2}} \hat{h}_{2,k}^{(3)} \cdot P_k^{(3 \times 3)} \cdot \hat{h}_{2,k}^{(3)}\right)$$
with a constant $\mu>0$ from the one species case done in \cite{achleitner2017multi}.

The cases $M>3$ are analogue to the case $M=3$ since then we have the same structure as in \eqref{eq10}, only with more entries $1$ on the diagonal in the term coming from the right-hand side of \eqref{BGKend}. It also reduces to the one species case.
\subsubsection{The case k=0}
For $k=0$, the term corresponding to the $x$-derivative vanishes, and the remaining terms with $P_0= \textbf{1}$ simplify to
\begin{align}
\begin{split}
&\frac{d}{dt} \hat{h}_{1,0} \cdot \hat{h}_{1,0} + \hat{h}_{1,0} \cdot \frac{d}{dt} \hat{h}_{1,0} = - 2 \nu_{11} n_{\infty,1}[ |\hat{h}_{1,(0,3)}|^2 + |\hat{h}_{1,(0,4)}|^2+ \cdots ] - 2 \nu_{12} n_{\infty,2} (1- \delta) |\hat{h}_{1,(0,1)}|^2\\&- 2 \nu_{12} n_{\infty,2} (1- \alpha) |\hat{h}_{1,(0,2)}|^2 - 2 \nu_{12} n_{\infty,2} |\hat{h}_{1,(0,3)}|^2 + |\hat{h}_{1,(0,4)}|^2 + \cdots ]  \\&+ 2 \nu_{12} n_{\infty,2} (1- \delta) \frac{n_{\infty,1}}{n_{\infty,2}} \frac{\sqrt{m_1}}{\sqrt{m_2}} \text{Re}( \overline{\hat{h}_{1,(0,1)}} \hat{h}_{2,(0,1)}) + 2 \nu_{12} n_{\infty,2} \frac{n_{\infty,1}}{n_{\infty,2}} (1- \alpha) \text{Re}(\hat{h}_{2,(0,2)} \overline{\hat{h}_{1,(0,2)}} )
\end{split}
\label{gl1}
\end{align}
for species $1$, and 
\begin{align}
\begin{split}
&\frac{d}{dt} \hat{h}_{2,0} \cdot \hat{h}_{2,0} + \hat{h}_{2,0} \cdot \frac{d}{dt} \hat{h}_{2,0} = - 2 \nu_{22} n_{\infty,2} [ |\hat{h}_{2,(0,3)}|^2 + |\hat{h}_{2,(0,4)}|^2 + \cdots ] - 2 \nu_{12} n_{\infty,1} \frac{m_1}{m_2} (1- \delta) |\hat{h}_{2,(0,1)}|^2 \\[4pt]&- 2 \nu_{12} n_{\infty,1} (1- \alpha) |\hat{h}_{2,(0,2)}|^2  - 2 \nu_{12} n_{\infty,1} [ |\hat{h}_{2,(0,3)}|^2 + |\hat{h}_{2,(0,4)}|^2 + \cdots] \\&+ 2 \nu_{12} n_{\infty,1} (1- \delta) \frac{n_{\infty,2}}{n_{\infty,1}} \frac{\sqrt{m_1}}{\sqrt{m_2}} \text{Re}( \overline{\hat{h}_{1,(0,1)}} \hat{h}_{2,(0,1)}) + 2 \nu_{12} n_{\infty,1} \frac{n_{\infty,2}}{n_{\infty,1}} (1- \alpha) \text{Re}(\hat{h}_{2,(0,2)} \overline{\hat{h}_{1,(0,2)}}) 
\end{split}
\label{gl2}
\end{align}
for species $2$. 
Due to conservation of the number of particles \eqref{cons} for each species, we have
$$ \hat{h}_{1,(0,0)}= \hat{h}_{2,(0,0)}=0.$$
Therefore, we can add 
$$ - C~ \frac{1}{n_{\infty,1}} ~|\hat{h}_{1,(0,0)}|^2 - C~ \frac{1}{n_{\infty,2}} ~|\hat{h}_{2,(0,0)}|^2$$
for an arbitrary positive constant $C>0$ to the sum of $\frac{1}{n_{\infty,1}}$ \eqref{gl1} and $\frac{1}{n_{\infty,2}}$ \eqref{gl2}, then we get
\begin{align}
\begin{split}
&\quad\frac{1}{n_{\infty,1}} (\frac{d}{dt} \hat{h}_{1,0} \cdot \hat{h}_{1,0} + \hat{h}_{1,0} \cdot \frac{d}{dt} \hat{h}_{1,0})+\frac{1}{n_{\infty,2}}(\frac{d}{dt}\hat{h}_{2,0} \cdot \hat{h}_{2,0} + \hat{h}_{2,0} \cdot \frac{d}{dt} \hat{h}_{2,0}) 
\\[4pt]&= - 2 \nu_{11} [ |\hat{h}_{1,(0,3)}|^2 + |\hat{h}_{1,(0,4)}|^2+ \cdots ] - 2 \nu_{12} \frac{n_{\infty,2}}{n_{\infty,1}} (1- \delta) |\hat{h}_{1,(0,1)}|^2\\&- 2 \nu_{12} \frac{n_{\infty,2}}{n_{\infty,1}} (1- \alpha) |\hat{h}_{1,(0,2)}|^2 - 2 \nu_{12} \frac{n_{\infty,2}}{n_{\infty,1}} |\hat{h}_{1,(0,3)}|^2 + |\hat{h}_{1,(0,4)}|^2 + \cdots ]  \\&+ 2 \nu_{12}  (1- \delta)  \frac{\sqrt{m_1}}{\sqrt{m_2}} \text{Re}( \hat{h}_{1,(0,1)} \hat{h}_{2,(0,1)}) + 2 \nu_{12}  (1- \alpha) \text{Re}(\hat{h}_{2,(0,2)} \hat{h}_{1,(0,2)} )
\\[4pt]&\quad- 2 \nu_{22}  [ |\hat{h}_{2,(0,3)}|^2 + |\hat{h}_{2,(0,4)}|^2 + \cdots ] - 2 \nu_{12} \frac{n_{\infty,1}}{n_{\infty,2}} \frac{m_1}{m_2} (1- \delta) |\hat{h}_{2,(0,1)}|^2 \\[4pt]&- 2 \nu_{12} \frac{n_{\infty,1}}{n_{\infty,2}} (1- \alpha) |\hat{h}_{2,(0,2)}|^2  - 2 \nu_{12} \frac{n_{\infty,1}}{n_{\infty,2}} [ |\hat{h}_{2,(0,3)}|^2 + |\hat{h}_{2,(0,4)}|^2 + \cdots] \\&+ 2 \nu_{12}  (1- \delta)  \frac{\sqrt{m_1}}{\sqrt{m_2}} \text{Re}( \overline{\hat{h}_{1,(0,1)}} \hat{h}_{2,(0,1)}) + 2 \nu_{12}   (1- \alpha) \text{Re}(\hat{h}_{2,(0,2)} \overline{\hat{h}_{1,(0,2)}}) \\ &- C~ \frac{1}{n_{\infty,1}} ~|\hat{h}_{1,(0,0)}|^2 - C~ \frac{1}{n_{\infty,2}} ~|\hat{h}_{2,(0,0)}|^2
\end{split}
\label{gl3}
\end{align}
Due to conservation of total momentum and energy \eqref{cons}, we have
$$\sqrt{m_1} \hat{h}_{1,(0,1)} + \sqrt{m_2} \hat{h}_{2,(0,1)} =0, \quad \hat{h}_{1,(0,2)} + \hat{h}_{2,(0,2)} = 0.$$
Multiplying them by $\sqrt{m_1} \overline{\hat{h}_{1,(0,1)}} + \sqrt{m_2} \overline{\hat{h}_{2,(0,1)}}$ and $\overline{\hat{h}_{1,(0,2)}} + \overline{\hat{h}_{2,(0,2)}}$ respectively, 
we observe that $\text{Re}(\overline{\hat{h}_{1,(0,1)}} \hat{h}_{2,(0,1)})$ and $\text{Re}(\overline{\hat{h}_{1,(0,2)}} \hat{h}_{2,(0,2)})$ must be negative. Therefore, we can estimate the right-hand side of \eqref{gl3} from above by
\begin{align}
\begin{split}
&- 2 \nu_{11} [ |\hat{h}_{1,(0,3)}|^2 + \cdots] - 2 \nu_{12}  \frac{n_{\infty,2}}{n_{\infty,1}} (1- \delta) |\hat{h}_{1,(0,1)}|^2- 2 \nu_{12}  \frac{n_{\infty,2}}{n_{\infty,1}} (1- \alpha) |\hat{h}_{1,(0,2)}|^2 - 2 \nu_{12}  \frac{n_{\infty,2}}{n_{\infty,1}} [|\hat{h}_{1,(0,3)}|^2 + \cdots] \\[2pt] &- 2 \nu_{22} [ \hat{h}^2_{2,(0,3)} + \cdots ] - 2 \nu_{12}  \frac{n_{\infty,1}}{n_{\infty,2}} \frac{m_1}{m_2} (1- \delta) |\hat{h}_{2,(0,1)}|^2- 2 \nu_{12}  \frac{n_{\infty,1}}{n_{\infty,2}} (1- \alpha) |\hat{h}_{2,(0,2)}|^2 \\[2pt] &- 2 \nu_{12}  \frac{n_{\infty,1}}{n_{\infty,2}} [|\hat{h}_{2,(0,3)}|^2 + \cdots]  - C~\frac{1}{n_{\infty,1}}~ |\hat{h}_{1,(0,0)}|^2 - C~ \frac{1}{n_{\infty,2}}~|\hat{h}_{2,(0,0)}|^2.
\label{gl4}
\end{split}
\end{align}
We choose {\footnotesize $$C= 2 \min \lbrace \nu_{12} n_{\infty,2} (1- \delta), \nu_{12} n_{\infty,2} (1-\alpha), \nu_{11} n_{\infty,1}+ \nu_{12} n_{\infty,2}, \nu_{12} n_{\infty,1} \frac{m_1}{m_2} (1- \delta), \nu_{12} n_{\infty,1} (1- \alpha), \nu_{22} n_{\infty,2} + \nu_{12} n_{\infty,1} \rbrace.$$}

In addition, we estimate all coefficients in \eqref{gl4} from below by $C$. Then, we obtain
\begin{align}
\begin{split}
&\quad\frac{1}{n_{\infty,1}} \left(\frac{d}{dt} \hat{h}_{1,0} \cdot \hat{h}_{1,0} + \hat{h}_{1,0} \cdot \frac{d}{dt} \hat{h}_{1,0}\right)+\frac{1}{n_{\infty,2}}\left(\frac{d}{dt} \hat{h}_{2,0} \cdot \hat{h}_{2,0} + \hat{h}_{2,0} \cdot \frac{d}{dt} \hat{h}_{2,0}\right)\\[2pt]& \leq - C \left(\frac{1}{n_{\infty,1}}(\hat{h}^2_{1,(0,0)}+ \hat{h}^2_{1,(0,1)}+ \hat{h}^2_{1,(0,2)} + \hat{h}^2_{1,(0,3)})+\frac{1}{n_{\infty,2}} (\hat{h}^2_{2,(0,0)}  + \hat{h}^2_{2,(0,1)}  + \hat{h}^2_{2,(0,2)}  + \hat{h}^2_{2,(0,3)}) \right)
\end{split}
\end{align}
 
In conclusion, we obtain that 
$$\frac{1}{n_{\infty,1}}\left(\frac{d}{dt}\hat{h}_{1,k} \cdot P_k \cdot \hat{h}_{1,k}+\hat{h}_{1,k} \cdot P_k \cdot \frac{d}{dt}\hat{h}_{1,k}\right) + \frac{1}{n_{\infty,1}}\left(\frac{d}{dt} \hat{h}_{2,k} \cdot P_k \cdot  \hat{h}_{2,k}+ \hat{h}_{2,k} \cdot P_k \cdot \frac{d}{dt} \hat{h}_{2,k}\right)$$
can be estimated from above by
$$- C \left(\frac{1}{n_{\infty,1}}(\hat{h}_{1,0} \cdot P_0  \cdot \hat{h}_{1,0}) + \frac{1}{n_{\infty,2}}(\hat{h}_{2,0} \cdot P_0 \cdot \hat{h}_{2,0})\right)$$
for $k=0$.

\section{Proof of Theorem \ref{theo}}
The first statement of theorem \ref{theo} is basically proven in \cite{achleitner2017multi}. We just take a linear combination of the two entropies of the two species.

It remains to prove the second statement of theorem \ref{theo}.
In the previous section we proved that for a fixed $M\geq3$, one has
\begin{align*}
&\quad\frac{d}{dt}\left(\frac{1}{n_{\infty,1}}(\hat{h}_{1,k}^{(M)} \cdot P_k^{(M \times M)} \cdot \hat{h}_{1,k}^{(M)})+ \frac{1}{n_{\infty,2}}(\hat{h}_{2,k}^{(M)} \cdot P_k^{(M \times M)} \cdot \hat{h}_{2,k}^{(M)})\right)\\[4pt]&=\frac{1}{n_{\infty,1}}\left(\frac{d}{dt}\hat{h}_{1,k}^{(M)} \cdot P_k^{(M \times M)} \cdot \hat{h}_{1,k}^{(M)}+\hat{h}_{1,k}^{(M)} \cdot P_k^{(M \times M)} \cdot \frac{d}{dt}\hat{h}_{1,k}^{(M)}\right) \\[4pt]&\quad + \frac{1}{n_{\infty,2}}\left(\frac{d}{dt} \hat{h}_{2,k}^{(M)} \cdot P_k^{(M \times M)} \cdot  \hat{h}_{2,k}^{(M)}+ \hat{h}_{2,k}^{(M)} \cdot P_k^{(M \times M)} \cdot \frac{d}{dt} \hat{h}_{2,k}^{(M)}\right)
\\[4pt] &\leq - 2 \mu \left(\frac{1}{n_{\infty,1}}(\hat{h}_{1,k}^{(M)} \cdot P_k^{(M \times M)} \cdot \hat{h}_{1,k}^{(M)}) + \frac{1}{n_{\infty,1}}(\hat{h}_{2,k}^{(M)} \cdot P_k^{(M \times M)} \cdot \hat{h}_{2,k}^{(M)})\right)
\end{align*}
for $k>0$, and
 \begin{align*}
 &\quad\frac{d}{dt}\left(\frac{1}{n_{\infty,1}}(\hat{h}_{1,0} \cdot P_0 \cdot \hat{h}_{1,0})+ \frac{1}{n_{\infty,2}}(\hat{h}_{2,0} \cdot P_0 \cdot \hat{h}_{2,0})\right) \\[4pt] &=
\frac{1}{n_{\infty,1}}\left(\frac{d}{dt}\hat{h}_{1,0} \cdot P_0\cdot \hat{h}_{1,0}+\hat{h}_{1,0} \cdot P_0 \cdot \frac{d}{dt}\hat{h}_{1,0}\right) \\[4pt]&\quad + \frac{1}{n_{\infty,2}}\left(\frac{d}{dt} \hat{h}_{2,0} \cdot P_0\cdot  \hat{h}_{2,0}+ \hat{h}_{2,0} \cdot P_0 \cdot \frac{d}{dt} \hat{h}_{2,0}\right) \\[4pt] &\leq - C \left(\frac{1}{n_{\infty,1}}(\hat{h}_{1,0}\cdot P_0 \cdot \hat{h}_{1,0}) + \frac{1}{n_{\infty,2}}(\hat{h}_{2,0} \cdot P_0 \cdot \hat{h}_{2,0}) \right) 
\end{align*}
for $k=0$. So we can deduce with Gronwall's lemma
\begin{align*}
&\quad\frac{1}{n_{\infty,1}}\left(\hat{h}_{1,k}^{(M)} \cdot P_k^{(M \times M)} \cdot \hat{h}_{1,k}^{(M)}\right)+ \frac{1}{n_{\infty,2}}\left(\hat{h}_{2,k}^{(M)} \cdot P_k^{(M \times M)} \cdot \hat{h}_{2,k}^{(M)}\right)\\[4pt] & \leq e^{- 2 \mu t} \left(\frac{1}{n_{\infty,1}}(\hat{h}_{1,k}^{(M)}(0) \cdot P_k^{(M \times M)} \cdot \hat{h}_{1,k}^{(M)}(0)) + \frac{1}{n_{\infty,2}}(\hat{h}_{2,k}^{(M)}(0) \cdot P_k^{(M \times M)} \cdot \hat{h}_{2,k}^{(M)}(0))\right) \\[4pt]
\end{align*}
for $k>0$, and
\begin{align*}
&\quad\frac{1}{n_{\infty,1}}\left(\hat{h}_{1,0} \cdot P_0 \cdot \hat{h}_{1,0}\right)
+ \frac{1}{n_{\infty,2}}\left(\hat{h}_{2,0} \cdot P_0 \cdot \hat{h}_{2,0}\right)\\[4pt]
&\leq e^{- C t}\left(\frac{1}{n_{\infty,1}}(\hat{h}_{1,0}(0) \cdot P_0\cdot \hat{h}_{1,0}(0)) + \frac{1}{n_{\infty,2}}(\hat{h}_{2,0}(0) \cdot P_0 \cdot \hat{h}_{2,0}(0))\right) 
\end{align*}
for $k=0$. Therefore, we have
\begin{align*}
&\frac{1}{n_{\infty,1}} e_{k,1}(\tilde{f_1}) + \frac{1}{n_{\infty,2}} e_{k,2}(\tilde{f}_2) \\&\leq e^{-2 \mu t} \lim_{M \rightarrow \infty} \left(\frac{1}{n_{\infty,1}}(\hat{h}_{1,k}^{(M)}(0) \cdot P_k^{(M \times M)} \cdot \hat{h}_{1,k}^{(M)}(0)) + \frac{1}{n_{\infty,2}}(\hat{h}_{2,k}^{(M)}(0) \cdot P_k^{(M \times M)} \cdot \hat{h}_{2,k}^{(M)}(0))\right)\\[4pt] 
& \leq e^{-2 \mu t}  \left( \frac{1}{n_{\infty,1}}(\hat{h}_{1,k}(0) \cdot P_k \cdot \hat{h}_{1,k}(0)) + \frac{1}{n_{\infty,2}}(\hat{h}_{2,k}(0) \cdot P_k \cdot \hat{h}_{2,k}(0))\right)
\end{align*}
according to the identity of Bessel for $k>0$, and
\begin{align*}
\frac{1}{n_{\infty,1}} e_{0,1}(\tilde{f_1}) + \frac{1}{n_{\infty,2}} e_{0,2}(\tilde{f}_2) & \leq e^{- C t} \left(\frac{1}{n_{\infty,1}}e_{0,1}(\tilde{f}_1(0)) + \frac{1}{n_{\infty,2}}e_{0,2}(\tilde{f}_2(0))\right)
\end{align*}
for $k=0$. Finally, this leads to
\begin{align*}
& e(\tilde{f}_1, \tilde{f}_2) = \sum_{k \in \mathbb{Z}} (\frac{1}{n_{\infty,1}} e_{k,1}(\tilde{f}_1) +\frac{1}{n_{\infty,2}} e_{k,2} (\tilde{f}_2)) \\& 
\qquad\quad\leq e^{- \min \lbrace C, 2 \mu \rbrace t}   e(\tilde{f}_1(0), \tilde{f}_2(0)). 
\end{align*}

Thus, we proved Theorem \ref{theo}.

\section{Conclusion and future work}
We considered a kinetic model for a two component gas mixture without chemical reactions. We studied hypocoercivity for the linearized BGK model for gas mixtures in continuous phase space. By constructing an entropy functional, 
we could prove exponential relaxation to equilibrium with explicit rates for the mixture system. 

As for the future work, we propose to extend the hypocoercivity estimates in the current work to the random case, that is, to conduct sensitivity analysis for the mixture BGK model with random inputs. Uncertainties may come from the initial data, various parameters in the model. To numerically solve the mixture BGK model with random inputs, a generalized polynomial chaos based stochastic Galerkin (gPC-SG) method can be used \cite{XiuBook}. It would be interesting to obtain estimates for the underlying gPC-SG system, and study spectral accuracy and exponential decay in time of the numerical error of the gPC method. Similar analysis for a general class of collisional kinetic equations with multiple scales and random inputs has been studied in \cite{LJ-UQ}. 

\section*{Acknowledgements}
Both authors would like to thank Prof. Christian Klingenberg and Prof. Shi Jin for their discussions and bringing authors
interests to work on this topic. Liu Liu is supported by the funding DOE--Simulation Center for Runaway Electron Avoidance and Mitigation, project number DE-SC0016283. Marlies Pirner is supported by the Austrian Science Fund (FWF) project F 65 and the Humboldt foundation. 

\section*{References:}
\bibliographystyle{siam}
\bibliography{multispecies.bib}

\end{document}